\documentclass[twoside,11pt]{article}

\usepackage[preprint]{jmlr2e}

\usepackage{jmlr2e}
\usepackage{amsfonts}
\usepackage{amsmath}
\usepackage{color}
\definecolor{bred}{rgb}{0.8,0,0}

\usepackage{hyperref}
\hypersetup{colorlinks,linkcolor={blue},citecolor={bred},urlcolor={blue}}

\def\z{{\mathbf z}}

\def\cB{{\mathcal B}}

\def\cF{{\mathcal F}}
\def\cG{{\mathcal G}}

\def\cC{{\mathcal C}}
\def\cV{{\mathcal V}}

\def\bR{{\mathbb R}}
\def\bU{{\mathbb U}}

\def\bE{{\mathbb E}}
\def\E{{\mathbb E}}

\def\cP{{\mathcal P}}
\def\cB{{\mathcal B}}
\def\bN{{\mathbb N}}

\def\cL{{\mathcal L}}
\def\f0{{\mathbf 0}}

\def\md{{\mathrm d}}

\def\sfp{{\sf p}}

\newcommand{\floor}[1]{{\lfloor #1 \rfloor}}
\newcommand{\ceil}[1]{{\lceil #1 \rceil}}

\newtheorem{thm}{Theorem}[section]

\newtheorem{defn}{Definition}[section]
\newtheorem{assumption}{Assumption}[section]

\newtheorem{lem}{Lemma}[section]
\newtheorem{cor}{Corollary}[section]
\newtheorem{prop}{Proposition}[section]

\newtheorem{rem}{Remark}[section]

\jmlrheading{1}{2024}{1-48}{4/00}{10/00}{akyildiz23a}{\"Omer Deniz Akyildiz and Sotirios Sabanis}

\ShortHeadings{Nonasymptotic analysis of SGHMC}{Akyildiz and Sabanis}
\firstpageno{1}

\begin{document}

\title{Nonasymptotic analysis of Stochastic Gradient Hamiltonian Monte Carlo under local conditions for nonconvex optimization}

\author{\name \"O. Deniz Akyildiz \email deniz.akyildiz@imperial.ac.uk \\
       \addr Department of Mathematics, Imperial College London, UK
       \AND
       \name Sotirios Sabanis \email s.sabanis@ed.ac.uk \\
       \addr School of Mathematics, University of Edinburgh, UK \\ The Alan Turing Institute, UK \\ National Technical University of Athens, Greece}


\maketitle

\begin{abstract}
We provide a nonasymptotic analysis of the convergence of the stochastic gradient Hamiltonian Monte Carlo (SGHMC) to a target measure in Wasserstein-2 distance without assuming log-concavity. Our analysis quantifies key theoretical properties of the SGHMC as a sampler under local conditions which significantly improves the findings of previous results. In particular, we prove that the Wasserstein-2 distance between the target and the law of the SGHMC is uniformly controlled by the step-size of the algorithm, therefore demonstrate that the SGHMC can provide high-precision results uniformly in the number of iterations. The analysis also allows us to obtain nonasymptotic bounds for nonconvex optimization problems under local conditions and implies that the SGHMC, when viewed as a nonconvex optimizer, converges to a global minimum with the best known rates. We apply our results to obtain nonasymptotic bounds for scalable Bayesian inference and nonasymptotic generalization bounds.
\end{abstract}


\section{Introduction}
\label{intro}
We are interested in nonasymptotic estimates for the sampling problem from the probability measures of the form
\begin{align}\label{eq:TargetMeasure}
\pi_\beta(\md\theta) \propto \exp(-\beta U(\theta)) \md \theta.
\end{align}
when only the noisy estimate of $\nabla U$ is available. This problem arises in many cases in machine learning, most notably in large-scale (mini-batch) Bayesian inference \citep{welling2011bayesian,ahn2012bayesian} and nonconvex stochastic optimization as for large values of $\beta$, a sample from the target measure \eqref{eq:TargetMeasure} is an approximate minimizer of the potential $U$ \citep{gelfand1991recursive}. Consequently, nonasymptotic error bounds for the sampling schemes can be used to obtain guarantees for Bayesian inference or nonconvex optimization.

An efficient method for obtaining a sample from \eqref{eq:TargetMeasure} is simulating the overdamped Langevin stochastic differential equation (SDE) which is given by
\begin{align}\label{eq:OverDampLangevinSDE}
\md L_t = -h(L_t) \md t + \sqrt{\frac{2}{\beta}} \md B_t,
\end{align}
with a random initial condition $L_0 := \theta_0$ where $h := \nabla U$ and $(B_t)_{t\geq 0}$ is a $d$-dimensional Brownian motion. The Langevin SDE \eqref{eq:OverDampLangevinSDE} admits $\pi_\beta$ as the unique invariant measure, therefore simulating this process will lead to samples from $\pi_\beta$ and can be used as a Markov chain Monte Carlo (MCMC) algorithm \citep{roberts1996exponential,roberts2002langevin}. Moreover, the fact that the limiting probability measure $\pi_\beta$ concentrates around the global minimum of $U$ for sufficiently large values of $\beta$ makes the diffusion \eqref{eq:OverDampLangevinSDE} also an attractive candidate as a global optimizer (see, e.g., \citet{hwang1980laplace}). However, since the continuous-time process \eqref{eq:OverDampLangevinSDE} cannot be simulated, its first-order Euler discretization with the step-size $\eta>0$ is used in practice, termed \textit{the Unadjusted Langevin Algorithm} (ULA) \citep{roberts1996exponential}. The ULA scheme has become popular in recent years due to its advantages in high-dimensional settings and ease of implementation. Nonasymptotic properties of the ULA were recently established under strong convexity and smoothness assumptions by \cite{dalalyan2017theoretical,durmus2017nonasymptotic, durmus2019high} while some extensions about relaxing smoothness assumptions or inaccurate gradients were also considered by \cite{dalalyan2019user,brosse2019tamed}. The similar attractive properties hold for the ULA when the potential $U$ is nonconvex \citep{gelfand1991recursive,raginsky2017non,xu2018global,erdogdu2018global, hola}. In recent works, ULA has been extended for nonconvex cases in several different directions, e.g., under log-Sobolev inequality \citep{vempala2019rapid}, under H\"older continuity and specific tail growth conditions \citep{erdogdu2021convergence}, under Poincar\'e inequality \citep{chewi2022analysis}. Further work has extended these results, see, e.g., \citet{balasubramanian2022towards} for averaged Langevin schemes, \citet{erdogdu2022convergence} for the analysis under dissipativity in Chi-squared and Renyi divergences, \citet{mou2022improved} for results under dissipativity with smooth initalisation, and finally under weak Poincar\'e inequalities \citep{mousavi2023towards}.

While the ULA performs well when the computation of the gradient $h(\cdot)$ is straightforward, this is not the case in most interesting applications. Usually, a stochastic, unbiased estimate of $h(\cdot)$ is available, either because the cost function is defined as an expectation or as a finite sum. Using stochastic gradients in the ULA leads to another scheme called stochastic gradient Langevin dynamics (SGLD) \citep{welling2011bayesian}. The SGLD has been particularly popular in the fields of (i) \textit{large-scale Bayesian inference} since it allows one to construct Markov chains Monte Carlo (MCMC) algorithms using only subsets of the dataset \citep{welling2011bayesian}, (ii) \textit{nonconvex optimization} since it enables one to estimate global minima using only stochastic (often cheap-to-compute) gradients \citep{raginsky2017non}. As a result, attempts for theoretical understanding of the SGLD have been recently made in several works, both for the strongly convex potentials (i.e. log-concave targets), see, e.g., \cite{barkhagen2021stochastic,brosse2018promises} and nonconvex potentials, see, e.g. \cite{raginsky2017non,majka2018non,zhang2023nonasymptotic}. Our particular interest is in nonasymptotic bounds for the nonconvex case, as it is relevant to our work. In their seminal paper, \cite{raginsky2017non} obtain a nonasymptotic bound between the law of the SGLD and the target measure in Wasserstein-2 distance with a rate $\eta^{5/4} n$ where $\eta$ is the step-size and $n$ is the number of iterations. While this work is first of its kind, the error rate grows with the number of iterations. In a related contribution, \cite{xu2018global} have obtained improved rates, albeit still growing with the number of iterations $n$. In a more recent work, \cite{chau2021stochastic} have obtained a uniform rate of order $\eta^{1/2}$ in Wasserstein-1 distance. \cite{majka2018non} achieved error rates of $\eta^{1/2}$ and $\eta^{1/4}$ for Wasserstein-1 and Wasserstein-2 distances, respectively, under the assumption of convexity outside a ball. Finally, \cite{zhang2023nonasymptotic} achieved the same rates under only local conditions which can be verified for a class of practical problems.

An alternative to methods based on the \textit{overdamped} Langevin SDE \eqref{eq:OverDampLangevinSDE} is the class of algorithms based on the \textit{underdamped} Langevin SDE. To be precise, the underdamped Langevin SDE is given as
\begin{align}
\md V_t &= -\gamma V_t \md t - h(\theta_t) \md t + \sqrt{\frac{2\gamma}{\beta}} \md B_t,\label{eq:UnderDampLangevinSDE1} \\
\md \theta_t &= V_t \md t \label{eq:UnderDampLangevinSDE2},
\end{align}
where $(\theta_t,V_t)_{t\geq 0}$ are called position and momentum process, respectively, and $h := \nabla U$. Similar to eq.~\eqref{eq:OverDampLangevinSDE}, this diffusion can be used as both an MCMC sampler and nonconvex optimizer, since under appropriate conditions, the Markov process $(\theta_t,V_t)_{t \geq 0}$ has a unique invariant measure given by
\begin{align}\label{eq:ExtendedTargetHMC}
\overline{\pi}_\beta(\md \theta, \md v) \propto \exp\left( -\beta \left( \frac{1}{2} \|v\|^2 + U(\theta) \right) \right) \md \theta \md v.
\end{align}
Consequently, the marginal distribution of \eqref{eq:ExtendedTargetHMC} in $\theta$ is precisely the target measure defined in \eqref{eq:TargetMeasure}. This means that sampling from \eqref{eq:ExtendedTargetHMC} in the extended space and then keeping the samples in the $\theta$-space would define a valid sampler for the sampling problem of \eqref{eq:TargetMeasure}.

Due to its attractive properties, methods based on the underdamped Langevin SDE have attracted significant attention. In particular, the first order discretization of \eqref{eq:UnderDampLangevinSDE1}--\eqref{eq:UnderDampLangevinSDE2}, which is termed \textit{underdamped Langevin MCMC} (i.e. the \textit{underdamped} counterpart of the ULA), has been a focus of attention, see, e.g., \cite{duncan2017using,dalalyan2018sampling,cheng2018underdamped}. Particularly, the underdamped Langevin MCMC has displayed improved convergence rates in the setting where $U$ is convex, see, e.g., \cite{dalalyan2018sampling,cheng2018underdamped}. Similar results have been extended to the nonconvex case. In particular, \cite{cheng2018sharp} have shown that the underdamped Langevin MCMC converges in Wasserstein-2 with a better dimension and step-size dependence under the assumptions smoothness and convexity outside a ball. It has been also shown that the underdamped Langevin MCMC can be seen as an accelerated optimization method in the space of measures in Kullback-Leibler divergence \citep{ma2019there}.

Similar to the case in the ULA, oftentimes $\nabla U(\cdot)$ is expensive or impossible to compute exactly, but rather an unbiased estimate of it can be obtained efficiently. The underdamped Langevin MCMC with stochastic gradients is dubbed as Stochastic Gradient Hamiltonian Monte Carlo (SGHMC) and given as \citep{chen2014stochastic,ma2015complete}
\begin{align}
{V}_{n+1}^\eta &= {V}_n^\eta - \eta [\gamma {V}_n^\eta + H({\theta}_n^\eta, X_{n+1})] + \sqrt{\frac{2\gamma\eta}{\beta}} \xi_{n+1}, \label{eq:SGHMC1} \\
{\theta}_{n+1}^\eta &= {\theta}_n^\eta + \eta {V}_n^\eta, \label{eq:SGHMC2}
\end{align}
where $\eta > 0$ is a step-size,  ${V}_0^\eta = v_0$, $ {\theta}_0^\eta = \theta_0$,  and $\bE[H(\theta,X_0)] = h(\theta)$ for every $\theta \in \bR^d$. We note that SGHMC relies on the Euler-Maruyama discretisation which is of the possible discretisation methods that can be used. Alternative discretisation methods can be considered and lead to different algorithms, e.g., see \citet{cheng2018underdamped, li2019stochastic, zhang2023improved}.

In this paper, we analyze recursions \eqref{eq:SGHMC1}--\eqref{eq:SGHMC2}. We achieve convergence bounds and improve the ones proved in \cite{gao2018global} and \cite{chau2022stochastic} (see Section~\ref{sec:Related} for a direct comparison).

\textbf{Notation.} For an integer $d\geq 1$, the Borel sigma-algebra of $\mathbb{R}^d$ is denoted by $\cB(\bR^d)$. We denote the dot product with $\langle \cdot,\cdot\rangle$ while $|\cdot|$ denotes the associated norm. The set of probability measures defined on a measurable space $(\bR^d,\mathcal{B}(\bR^d))$ is denoted as $\mathcal{P}(\bR^d)$. For an $\bR^d$-valued random variable, $\cL(X)$ and $\bE[X]$ are used to denote its law and its expectation respectively. Note that we also write $\bE[X]$  as $\bE X$. For $\mu,\nu \in \cP(\bR^d)$, let $\cC(\mu,\nu)$ denote the set of probability measures $\Gamma$ on $\cB(\bR^{2d})$ so that its marginals are $\mu,\nu$. Finally, $\mu,\nu\in\mathcal{P}(\bR^d)$, the Wasserstein distance of order $p\geq 1$ is defined as
\begin{equation}\label{eq:definition-W-p}
{W}_p(\mu,\nu):=\inf_{\Gamma\in\mathcal{C}(\mu,\nu)}
\left(\int_{\bR^d}\int_{\bR^d}|\theta-\theta'|^p\Gamma(\md \theta , \md \theta')\right)^{1/p}.
\end{equation}

\section{Main results and overview}
Let $(X_n)_{n\in\bN}$ be an $\bR^m$-valued stochastic process adapted to $(\cG_n)_{n\in\bN}$ where $\cG_n := \sigma(X_k, k\leq n, k\in \bN)$ for $n \in \bN$. It is assumed henceforth that $\theta_0,v_0, \cG_\infty$, and $(\xi_n)_{n\in\bN}$ are independent. The main assumptions about $U$ and other quantities follow.
\begin{assumption}\label{ass:nonnegativity}
The cost function $U$ takes nonnegative values, i.e., $U(\theta) \geq 0$.
\end{assumption}
Next, the \textit{local smoothness} assumptions on the stochastic gradients $H(\theta,\cdot)$ are given.
\begin{assumption}\label{loclip}
There exist positive constants $L_1$, $L_2$ and $\rho$ such that, for all $x,x'\in\mathbb{R}^m$ and $\theta, \theta'\in\mathbb{R}^d$,
\begin{align*}
  |H(\theta,x)- H(\theta',x)| & \le L_1(1+|x|)^{\rho}|\theta-\theta'|, \\
  |H(\theta,x)- H(\theta,x')| & \le L_2(1+|x|+|x'|)^{\rho}(1+ |\theta|)|x-x'|.
\end{align*}
\end{assumption}
The following assumption states that the stochastic gradients are assumed to be unbiased.
\begin{assumption}\label{assmp:iid}
The process $(X_n)_{n \in \bN}$ is i.i.d. with $|X_0| \in L^{{4}(\rho +1)}$ and $|\theta_0|, |v_0| \in L^{{4}}$. It satisfies
\[
\bE[H(\theta,X_0)]=h(\theta).
\]
\end{assumption}

It is important to note that Assumption~\ref{loclip} is a significant relaxation in comparison with the corresponding assumptions provided in the literature, see, e.g., \cite{raginsky2017non,gao2018global,chau2022stochastic}. To the best of the authors' knowledge, all relevant works in this area have focused on uniform Lipschitz assumptions with the exception of \cite{zhang2023nonasymptotic}, which provides a nonasymptotic analysis of the SGLD under similar assumptions to ours.
\begin{rem}\label{rem:BoundsOnH} Assumption \ref{loclip} implies, for all $\theta,\theta' \in \bR^{d}$,
\begin{equation}\label{eq:LipschitzGradient}
| h(\theta)-h(\theta')| \leq L_1\bE[(1+ |X_0|)^{\rho}]|\theta-\theta' |,
\end{equation}
which consequently implies
\begin{align}\label{eq:BoundedGradh}
|h(\theta)| \leq L_1 \bE[(1 + |X_0|)^\rho] | \theta| + h_0,
\end{align}
where $h_0 := |h(0)|$. Let $H_0 := |H(0,0)|$, then Assumption \ref{loclip} implies
\[
|H(\theta,x)|\leq L_1(1+|x|)^{\rho}|\theta|+L_2(1+|x|)^{\rho+1}+H_0.
\]
\end{rem}
We denote $C_\rho := \bE \left[ (1 + |X_0|)^{{4}(\rho + 1)} \right]$. Note that $C_\rho < \infty$ by Assumption~\ref{assmp:iid}.
\begin{assumption}\label{assum:dissipativity}
{There exist a measurable (symmetric matrix-valued) function} $A:\bR^m\to\bR^{d\times d}$, $b: \mathbb{R}^m \to \mathbb{R}$ such that for any $ x,y \in \bR^d$, $\langle y, A(x) y\rangle \geq 0$, and for all $\theta \in \mathbb{R}^d$ and $x\in\mathbb{R}^m$,
\begin{align*}
\langle H(\theta,x),\theta\rangle\geq \langle \theta, A(x) \theta\rangle -b(x).
\end{align*}
The smallest eigenvalue of $\bE[A(X_0)]$ is a positive real number $a>0$ and $\bE[b(X_0)] = b>0$.
\end{assumption}
Note that Assumption~\ref{assum:dissipativity} is a \textit{local} dissipativity condition. This assumption implies the usual dissipativity property on the corresponding deterministic (full) gradient. In the next remark, we motivate these assumptions for the case of linear regression loss.
{\begin{rem}\label{rem:linear_regression} We note that Assumptions~\ref{loclip} and \ref{assum:dissipativity} can be motivated even using the simplest linear regression loss. Consider the problem
\begin{align*}
\min_{\theta \in \bR^d} \bE \left[ |Z - \langle Y, \theta \rangle |^2 \right],
\end{align*}
where $(Z, Y) \sim P(\mathrm{d} z, \md y)$ on $\bR \times \bR^d$. The stochastic gradient here is given by, for a sample $(z_n, y_n) \sim P(\md z, \md y)$,
\begin{align*}
H(\theta, x_n) = -2 y_n z_n - y_n \langle y_n, \theta \rangle,
\end{align*}
where $x_n = (z_n, y_n)$. This loss is \textit{not} globally Lipschitz however we can \textit{prove} that it is locally Lipschitz. For example, it is easy to show that Assumption~\ref{loclip} is satisfied with $L_1 = 1$, $L_2 = 2$, and $\rho = 2$. Moreover, it is easy to see that this stochastic gradient is not dissipative uniformly in $x_n$ as
\begin{align*}
\langle \theta, H(\theta, x_n) \rangle \geq |\langle y_n, \theta \rangle |^2 - z_n^2.
\end{align*}
However, our Assumption~\ref{assum:dissipativity} is satisfied with $A(x_n) = y_n y_n^\top$ and $b(x_n) = z_n^2$. We would also like to note that we do not need to assume that stochastic gradient moments are bounded as it is a result of our assumptions.
\end{rem}}
\begin{rem} \label{R2}
By Assumption \ref{assum:dissipativity}, we obtain $\left\langle h(\theta),\theta\right\rangle\geq a |\theta|^2-b$, for $\theta\in\mathbb{R}^d$ and $a,b > 0$.
\end{rem}
Below, we state our main result about the convergence of the law $\cL(\theta_k^\eta,V_k^\eta)$, which is generated by the SGHMC recursions \eqref{eq:SGHMC1}--\eqref{eq:SGHMC2}, to the extended target measure $\overline{\pi}_\beta$ in Wasserstein-2 ($W_2$) distance. We first define
\begin{align*}
\eta_{\textnormal{max}} = \min\left\lbrace 1, {\frac{1}{\gamma}}, \frac{\gamma\lambda}{2 K_1}, \frac{K_3}{K_2}, \frac{\gamma\lambda}{2\bar{K}}\right\rbrace,
\end{align*} 
where $K_1,K_2,K_3$ are constants explicitly given in the proof of Lemma~\ref{lem:DiscreteTimeLemma} and $\bar{K}$ is a constant explicitly given in the proof of Lemma~\ref{lem:SquareContraction}. Then, the following result is obtained.
\begin{thm}\label{thm:ConvRate} Let Assumptions \ref{ass:nonnegativity}--\ref{assum:dissipativity} hold. Then, there exist constants $C_1^\star, C_2^\star, C_3^\star, C_4^\star > 0$ such that, for every $0 < \eta \leq \eta_{\textnormal{max}}$,
\begin{align}
W_2(\cL(\theta_n^\eta,V_n^\eta),\overline{\pi}_\beta) &\leq C_1^\star {\eta}^{1/2} + C_2^\star \eta^{1/4} +C_3^\star e^{-C_4^\star \eta n}.
\end{align}
where the constants $C_1^\star,C_2^\star,C_3^\star,C_4^\star$ are explicitly provided in the Appendix.
\end{thm}
\begin{proof} See Section~\ref{sec:proof:thm:ConvRate}.
\end{proof}
\begin{rem}\label{rem:DimensionDependence}
We remark that $C_1^\star = \mathcal{O}(d^{1/2})$, $C_2^\star = \mathcal{O}(e^d)$, $C_3^\star = \mathcal{O}(e^d)$, and $C_4^\star = \mathcal{O}(e^{-d})$. We note that although the dependence of $C_1^\star$ on dimension is $\mathcal{O}(d^{1/2})$, and comes from our main result, the dependence of $C_2^\star$, $C_3^\star$, and $C_4^\star$ to dimension may be exponential as it is an immediate consequence of the contraction result of the underdamped Langevin SDE in \cite{eberle2019couplings}.
\end{rem}
The result in Theorem~\ref{thm:ConvRate} demonstrates that the error scales like $\mathcal{O}(\eta^{1/4})$ and is uniformly bounded over $n$ which can be made arbitrarily small by choosing $\eta > 0$ small enough. This result is thus a significant improvement over the findings in \cite{gao2018global}, where error bounds grow with the number of iterations, and in \cite{chau2022stochastic}, where the corresponding error bounds contain an additional term that is independent of $\eta$ and relates to the variance of the unbiased estimator.

\begin{rem} Our proof techniques can be adapted easily when $H(\theta,x) = h(\theta)$. Hence Theorem~\ref{thm:ConvRate} provides a novel convergence bound for the underdamped Langevin MCMC as well. {Our result is in $W_2$ distance under dissipativity for the Euler-Maruyama discretisation of the underdamped Langevin dynamics, which can be contrasted with existing work which focus on alternative discretisations. Beyond the log-concave case, \citet{ma2019there} analyses a second-order discretisation of underdamped Langevin diffusion under the assumption that the target satisfies a log-Sobolev inequality with a Lipschitz Hessian condition for the log-target. More recently, \citet{zhang2023improved} removed the Lipschitz Hessian condition and established the convergence of underdamped Langevin MCMC under log-Sobolev and Poincar\'e inequalities. To the best of our knowledge, our result for the Euler-Maruyama discretisation of the underdamped Langevin dynamics is the first nonasymptotic result under dissipativity that is uniformly controlled by the step-size.}
\end{rem}
Let $(\theta_k^\eta)_{k\in\bN}$ be generated by the SGHMC algorithm. Convergence of the $\cL(\theta_k^\eta)$ to $\pi_\beta$ in $W_2$ also implies that one can prove a global convergence result \citep{raginsky2017non}. More precisely, assume that we aim at solving the problem $\theta_\star \in \arg \min_{\theta\in\bR^d} U(\theta)$ which is a nonconvex optimization problem. We denote $U_\star := \inf_{\theta\in\bR^d} U(\theta)$.  Then we can bound the error $\bE[U(\theta_k^\eta)] - U_\star$ which would give us a guarantee on the nonconvex optimization problem. We state it as a formal result as follows.
\begin{thm}\label{thm:NonconvexBound} Under the assumptions of Theorem~\ref{thm:ConvRate}, we obtain
\begin{align*}
\bE[U(\theta_n^\eta)] - U_\star &\leq \overline{C}_1^\star {\eta}^{1/2} + \overline{C}_2^\star \eta^{1/4} +\overline{C}_3^\star e^{-C_4^\star \eta n} + \frac{d}{2\beta} \log \left( \frac{e \overline{L}_1}{a} \left(\frac{b\beta}{d} + 1\right) \right),
\end{align*}
where $\overline{C}_1^\star,\overline{C}_2^\star,\overline{C}_3^\star, C_4^\star, \overline{L}_1$ are finite constants which are explicitly given in the proofs.
\end{thm}
\begin{proof} See Section~\ref{sec:proof:thm:NonconvexBound}.
\end{proof}
This result bounds the error in terms of the function value for convergence to the global minima. Note that $\overline{C}_1^\star,\overline{C}_2^\star$ anf $\overline{C}_3^\star$ have the same dependence to dimension as $C_1^\star,C_2^\star$ and $C_3^\star$ (see Remark~\ref{rem:DimensionDependence}).

\subsection{Related work and contributions}\label{sec:Related}
Our work is related to two available analyses of the SGHMC, i.e.,  \cite{gao2018global} and \cite{chau2022stochastic}. We compare the bounds provided in Theorem~\ref{thm:ConvRate} and \ref{thm:NonconvexBound} to these two works. Finally, we also briefly consider the relations to \cite{zhang2023nonasymptotic} at the end of this section.

The scheme \eqref{eq:SGHMC1}--\eqref{eq:SGHMC2} is analyzed by \cite{gao2018global}. In particular, \cite{gao2018global} provided a convergence rate of the SGHMC \eqref{eq:SGHMC1}--\eqref{eq:SGHMC2} to the underdamped Langevin SDE \eqref{eq:UnderDampLangevinSDE1}--\eqref{eq:UnderDampLangevinSDE2} which is of order $\mathcal{O}(\delta^{1/4} + \eta^{1/4}) \sqrt{\eta n} \sqrt{\log(\eta n)}$. This rate grows with $n$, hence worsens over the number of iterations. Moreover, it is achieved under a uniform assumption on the stochastic gradient, i.e., $H(\theta,x)$ is assumed to be Lipschitz in $\theta$ uniformly in $x$ (as opposed to our Assumption~\ref{loclip}). Moreover, the mean-squared error of the gradient is assumed to be bounded whereas we do not place such an assumption in our work. Similar analyses appeared in the literature, e.g., for variance-reduced SGHMC \citep{zou2019stochastic} which also has growing rates with the number of iterations.

Another related work was provided by \cite{chau2022stochastic} who also analyzed the SGHMC recursions essentially under the same assumptions as in \cite{gao2018global}. However, \cite{chau2022stochastic} improved the convergence rate of the SGHMC recursions to the underdamped Langevin SDE significantly, i.e., provided a convergence rate of order $\mathcal{O}(\delta^{1/4} + \eta^{1/4})$ where $\delta > 0$ is a constant. While this rate significantly improves the rate of \cite{gao2018global}, it cannot be made to vanish by choosing $\eta > 0$ small enough, as $\delta > 0$ is (a priori assumed to be) independent of $\eta$.

In contrast, we prove that the SGHMC recursions track the underdamped Langevin SDE with a rate of order $\mathcal{O}(\eta^{1/4})$ which can be made arbitrarily small as with small $\eta > 0$. Moreover, our assumptions are significantly relaxed compared to \cite{gao2018global} and \cite{chau2022stochastic}. In particular, we relax the assumptions on stochastic gradients significantly by allowing growth in both variables $(\theta,x)$ which makes our theory hold for practical settings \citep{zhang2023nonasymptotic}.

Note that we analyze the SGHMC under similar assumptions to \cite{zhang2023nonasymptotic} who analyzed SGLD, since these assumptions allow for broader applicability of the results. However, the SGHMC requires the analysis of the underdamped Langevin diffusion and the discrete-time scheme on $\bR^{2d}$ and requires significantly different techniques compared to the ones used in \cite{zhang2023nonasymptotic}. We show that the SGHMC can achieve a similar performance with the SGLD as showed in \cite{zhang2023nonasymptotic}. However, due to the nonconvexity, we do not observe any improved dimension dependence as in the convex case and this constitutes a significant direction of future work.
\section{Preliminaries}\label{sec:Preliminaries}
In this section, preliminary results which are essential for proving the main results are provided. A central idea behind the proof of Theorem~\ref{thm:ConvRate} is the introduction of continuous-time auxiliary processes whose marginals at chosen discrete times coincide with the marginals of the (joint) law $\cL(\theta_k^\eta,V_k^\eta)$. Hence, these auxiliary stochastic processes can be used to analyze the recursions \eqref{eq:SGHMC1}--\eqref{eq:SGHMC2}. We first introduce the necessary lemmas about the moments of the auxiliary processes defined in Section~\ref{sec:Preliminaries} and contraction rates of the underlying underdamped Langevin diffusion. We give the lemmas and some explanatory remarks and defer the proofs to Appendix~\ref{app:sec:proof:prel_res}. We then provide the proofs of main theorems in the following section.

Consider the scaled process $({\zeta}_t^\eta,{Z}_t^{\eta}) := (\theta_{\eta t},V_{\eta t})$ given $(\theta_t,V_t)_{t\in\bR_+}$ as in \eqref{eq:UnderDampLangevinSDE1}--\eqref{eq:UnderDampLangevinSDE2}. We next define
\begin{align}
\md {Z}_t^{\eta} &= - \eta ( \gamma {Z}_t^\eta + h({\zeta}_t^\eta) ) \md t + \sqrt{2\gamma \eta \beta^{-1}} \md B_t^\eta, \\
\md {\zeta}_t^\eta &= \eta {Z}_t^\eta \md t,
\end{align}
where $\eta > 0$ and $B_t^\eta = \frac{1}{\sqrt{\eta}} B_{\eta t}$, where $(B_s)_{s\in\bR_+}$ is a Brownian motion with natural filtration $\cF_t$. We denote the natural filtration of $(B_t^\eta)_{t\in\bR_+}$ as $\cF_t^\eta$. We note that $\cF_t^\eta$ is independent of $\cG_\infty \vee \sigma(\theta_0,v_0)$. Next, we define the continuous-time interpolation of the SGHMC
\begin{align}
\md \overline{V}_t^\eta &= - \eta ( \gamma \overline{V}_\floor{t}^\eta + H(\overline{\theta}_\floor{t}^\eta , X_\ceil{t}) ) \md t + \sqrt{2\gamma \eta \beta^{-1}} \md B_t^\eta, \label{eq:ContTimeInterpSGHMC1} \\
\md \overline{\theta}_t^\eta &= \eta \overline{V}_\floor{t}^\eta \md t. \label{eq:ContTimeInterpSGHMC2}
\end{align}
The processes \eqref{eq:ContTimeInterpSGHMC1}--\eqref{eq:ContTimeInterpSGHMC2} mimic the recursions \eqref{eq:SGHMC1}--\eqref{eq:SGHMC2} at $n \in \bN$, i.e., $\cL(\theta_n^\eta,V_n^\eta) = \cL(\overline{\theta}_n^\eta,\overline{V}_n^\eta)$. Finally, we define the underdamped Langevin process $(\widehat{\zeta}_{t}^{s,u,v,\eta},\widehat{Z}_t^{s,u,v,\eta})$ for $s \leq t$
\begin{align}
\md \widehat{Z}_t^{s,u,v,\eta} &= - \eta ( \gamma \widehat{Z}_t^{s,u,v,\eta} + h(\widehat{\zeta}_{t}^{s,u,v,\eta}) ) \md t + \sqrt{2\gamma \eta \beta^{-1}} \md B_t^\eta, \\
\md \widehat{\zeta}_{t}^{s,u,v,\eta} &= \eta \widehat{Z}_t^{s,u,v,\eta} \md t,
\end{align}
with initial conditions $\widehat{\theta}_{s}^{s,u,v,\eta} = u$ and $\widehat{V}_s^{s,u,v,\eta} = v$.
\begin{defn} Fix $n\in\bN$ and for $T := \floor{1/\eta}$, define
\begin{align*}
\overline{\zeta}_t^{\eta,n} = \widehat{\zeta}_t^{nT,\overline{\theta}_{nT}^\eta,\overline{V}_{nT}^\eta,\eta}, \quad\quad \text{and} \quad \quad
\overline{Z}_t^{\eta,n} = \widehat{Z}_t^{nT,\overline{\theta}_{nT}^\eta,\overline{V}_{nT}\eta,\eta}.
\end{align*}
\end{defn}
The process $(\overline{\zeta}_t^{\eta,n},\overline{Z}_t^{\eta,n})_{t\ge nT}$ is an underdamped Langevin process started at time $nT$ with $(\overline{\theta}^\eta_{nT},\overline{V}^\eta_{nT})$.

To achieve the convergence results, we first define a Lyapunov function, common in the literature \citep{mattingly2002ergodicity,eberle2019couplings}, defined as
\begin{align}\label{eq:LyapunovDefn}
\mathcal{V}(\theta,v) = \beta U(\theta) + \frac{\beta}{4} \gamma^2 \left( |\theta + \gamma^{-1} v |^2 + |\gamma^{-1} v |^2 - \lambda |\theta |^2 \right),
\end{align}
where $\lambda \in (0,1/4]$. This Lyapunov function plays an important role in obtaining uniform moment estimates for some of the aforementioned processes. First we recall a result from the literature about the second moments of the processes $(\theta_t,V_t)_{t\geq 0}$ in continuous-time.
\begin{lem}\label{lem:ContTimeLemma} (Lemma~12(i) in \cite{gao2018global}.) Let Assumptions~\ref{ass:nonnegativity}--\ref{assum:dissipativity} hold. Then
\begin{align}
\sup_{t\geq 0} \bE | \theta_t |^2 \leq C_\theta^c &:= \frac{\int_{\bR^{2d}} \mathcal{V}(\theta,v) \mu_0(\md \theta,\md v) + \frac{d + A_c}{\lambda}}{\frac{1}{8} (1-2\lambda)\beta\gamma^2}, \\
\sup_{t\geq 0} \bE |V_t|^2 \leq C_v^c &:= \frac{\int_{\bR^{2d}} \mathcal{V}(\theta,v) \mu_0(\md \theta,\md v) + \frac{d + A_c}{\lambda}}{\frac{1}{4} (1-2\lambda)\beta}
\end{align}
\end{lem}
\begin{proof}
See \cite{gao2018global}.
\end{proof}
Next, we obtain, uniform in time, second moment estimates for the discrete-time processes $(\theta_k^\eta)_{k\geq 0}$ and $(V_k^\eta)_{k\geq 0}$.
\begin{lem}\label{lem:DiscreteTimeLemma} Let Assumptions~\ref{ass:nonnegativity}--\ref{assum:dissipativity} hold. Then, for $0 < \eta \leq \eta_\textnormal{max}$,
\begin{align*}
\sup_{k\geq 0} \bE |\theta_k^\eta|^2 \leq C_\theta &:= \frac{\int_{\bR^{2d}} \mathcal{V}(\theta,v) \mu_0(\md \theta, \md v) + \frac{4 (A_c + d)}{\lambda}}{\frac{1}{8}(1-2\lambda)\beta\gamma^2}, \\
\sup_{k \geq 0} \bE |V_k^\eta|^2 \leq C_v &:= \frac{\int_{\bR^{2d}} \mathcal{V}(\theta,v) \mu_0(\md \theta, \md v) + \frac{4 (A_c + d)}{\lambda}}{\frac{1}{4} (1- 2\lambda)\beta}.
\end{align*}
\end{lem}
\begin{proof}
See Appendix~\ref{proof:lem:DiscreteTimeLemma}.
\end{proof}
Using Lemma~\ref{lem:DiscreteTimeLemma}, we can get the second moment bounds of our auxiliary process $(\overline{\zeta}_t^{\eta,n})_{t\geq 0}$.
\begin{lem}\label{lem:zetaprocessmoment} Under the assumptions of Lemmas~\ref{lem:ContTimeLemma} and \ref{lem:DiscreteTimeLemma}, we obtain
\begin{align}
\sup_{n\in\bN} \sup_{t \in (nT, (n+1) T)} \bE | \overline{\zeta}_t^{\eta,n}|^2 \leq C_\zeta := \frac{\int_{\bR^{2d}} \mathcal{V}(\theta,v) \mu_0(\md \theta,\md v) + \frac{8 (d + A_c)}{\lambda}}{\frac{1}{8} (1-2\lambda)\beta\gamma^2}.
\end{align}
\end{lem}
\begin{proof}
See Appendix~\ref{proof:lem:zetaprocessmoment}.
\end{proof}
Moreover, in order to obtain the error bound distance between the laws $\cL(\overline{\zeta}_t^{\eta,n}, \overline{Z}_t^{\eta,n})$ and $\cL(\zeta_t^\eta, Z_t^\eta)$, we need to obtain the contraction of $\bE \cV^2(\theta_k^\eta, V_k^\eta)$, which is established in the following lemma.
\begin{lem}\label{lem:SquareContraction} Let $0 < \eta \leq \eta_{\textnormal{max}}$ and Assumptions~\ref{ass:nonnegativity}--\ref{assum:dissipativity} hold. Then, we have
\begin{align*}
\sup_{k\in\bN} \bE[\cV^2(\theta_k^\eta, V_k^\eta)] \leq \bE[\cV^2(\theta_0, v_0)] + \frac{2 D}{\gamma\lambda}
\end{align*}
where $D = \mathcal{O}(d^2)$ constant independent of $\eta$ and provided explicitly in the proof.
\end{lem}
\begin{proof}
See Appendix~\ref{proof:lem:SquareContraction}.
\end{proof}

Finally, we present a convergence result for the underdamped Langevin diffusion adapted from \cite{eberle2019couplings}. To this end, a functional for probability measures $\mu,\nu$ on $\bR^{2d}$ is introduced below
\begin{align}
\mathcal{W}_\rho(\mu,\nu) = \inf_{\Gamma\in\cC(\mu,\nu)} \int \rho((x,&v), (x',v')) \Gamma(d(x,v),d(x',v')),
\end{align}
where $\rho$ is defined in eq.~(2.10) in \cite{eberle2019couplings}. Thus, in view of Remarks \ref{rem:BoundsOnH} and \ref{R2}, one recovers the following result.
\begin{thm}\label{thm:Eberle} \citep[Theorem 2.3 and Corollary 2.6]{eberle2019couplings} Let Assumptions~\ref{ass:nonnegativity}--\ref{assum:dissipativity} hold and the laws of the underdamped Langevin SDEs $(\theta_t,V_t)$ and $(\theta'_t,V'_t)$ started at $(\theta_0,V_0) \sim \mu$ and $(\theta'_0,V'_0) \sim \nu$ respectively. Then, there exist constants $\dot{c},\dot{C} \in (0,\infty)$ such that
\begin{align}
W_2(\cL(\theta_t,V_t), \cL(\theta'_t,V'_t)) \leq \sqrt{\dot{C}} e^{-\dot{c}t / 2} \sqrt{\mathcal{W}_\rho(\mu,\nu)},
\end{align}
where the constants $\dot{c} = \mathcal{O}(e^{-d})$ and $\dot{C} = \mathcal{O}(e^d)$ are given in Appendix~\ref{app:Thm31}.
\end{thm}
\section{Proof of Theorem~\ref{thm:ConvRate}}\label{sec:proof:thm:ConvRate}
In order to prove Theorem~\ref{thm:ConvRate}, we note first that $W_2$ is a metric on $\mathcal{P}(\bR^{2d})$. The main strategy for this proof is to bound $W_2(\cL(\theta_n^\eta, V_n^\eta), \overline{\pi}_\beta)$ by using appropriate estimates on the continuous time interpolation of $(\theta_n^\eta,V_n^\eta)_{n\in\bN}$. In particular, we obtain the desired bound by decomposing first as
\begin{align}\label{eq:ErrDecomp}
W_2(\mathcal{L}(\overline{\theta}^\eta_t, \overline{V}_t^\eta), \overline{\pi}_\beta) &\leq W_2(\cL(\overline{\theta}_t^\eta,\overline{V}_t^\eta), \cL(\overline{\zeta}_t^{\eta,n}, \overline{Z}_t^{\eta,n})) + W_2(\cL(\overline{\zeta}_t^{\eta,n}, \overline{Z}_t^{\eta,n}), \cL(\zeta_t^\eta,Z_t^\eta)) \nonumber \\&+ W_2(\cL(\zeta_t^\eta,Z_t^\eta),\overline{\pi}_\beta),
\end{align}
for $nT \leq t \leq (n+1) T$ for every $n \in \bN$ and then by obtaining suitable (decaying in $n$) bounds for each of the terms on the rhs of \eqref{eq:ErrDecomp}. This leads to the proof of our main result, namely, Theorem~\ref{thm:ConvRate}. All proofs are deferred to the Appendix. First, we bound the first term of \eqref{eq:ErrDecomp}.
\begin{thm}\label{thm:W2_first_term} Let Assumptions~\ref{ass:nonnegativity}--\ref{assum:dissipativity} hold and $0<\eta\leq\eta_{\textnormal{max}}$. {Then, for every $t \in [nT, nT + 1)$,}
\begin{align}
W_2(\cL(\overline{\theta}^\eta_t,\overline{V}_t^\eta), \cL(\overline{\zeta}_t^{\eta,n}, \overline{Z}_t^{\eta,n})) \leq C_1^\star \eta^{1/2}
\end{align}
where $C_1^\star = \mathcal{O}(d^{1/2})$ is a finite constant.
\end{thm}
\begin{proof}
See Appendix~\ref{proof:thm:W2_first_term}.
\end{proof}
Next, we prove the following result for bounding the second term of \eqref{eq:ErrDecomp}.
\begin{thm}\label{thm:W2_second_term} Let Assumptions~\ref{ass:nonnegativity}--\ref{assum:dissipativity} hold and $0<\eta\leq\eta_{\textnormal{max}}$. Then,
\begin{align*}
W_2(\cL(\overline{\zeta}_t^{\eta,n}, \overline{Z}_t^{\eta,n}), \cL(\zeta_t^\eta,Z_t^\eta)) \leq C_2^\star \eta^{1/4},
\end{align*}
where $C_2^\star = \mathcal{O}(e^d)$.
\end{thm}
\begin{proof}
See Appendix~\ref{proof:thm:W2_second_term}.
\end{proof}
The constant $C_2^\star$ comes from the contraction result of \citet[Corollary~2.6]{eberle2019couplings} which might scale exponentially in $d$. Finally, the convergence of the last term follows from Theorem~\ref{thm:Eberle}.
\begin{thm}\label{thm:RateOfConvergenceScaledProcess} \citep{eberle2019couplings} Let Assumptions \ref{ass:nonnegativity}--\ref{assum:dissipativity} hold. Then,
\begin{align*}
W_2(\cL(\zeta_t^\eta,Z_t^\eta),\overline{\pi}_\beta) \leq C_3^\star e^{-C_4^\star \eta t},
\end{align*}
where $C_3^\star = \sqrt{\dot{C} \mathcal{W}_\rho(\mu_0,\nu_0)}$ and $C_4^\star = \dot{c}/2$. In particular, $C_3^\star = \mathcal{O}(e^d)$ while $C_4^\star = \mathcal{O}(e^{-d})$.
\end{thm}
Finally, considering Theorems~\ref{thm:W2_first_term}, \ref{thm:W2_second_term}, and \ref{thm:RateOfConvergenceScaledProcess} together by putting $t = n$ leads to the full proof of our main result, namely, Theorem~\ref{thm:ConvRate}.
\section{Proof of Theorem~\ref{thm:NonconvexBound}}\label{sec:proof:thm:NonconvexBound}
The bound provided for the convergence to the target in $W_2$ distance can be used to obtain guarantees for the nonconvex optimization. In order to do so, we proceed by decomposing the error as follows
\begin{align*}
\bE[U(\theta_n^\eta)] - U_\star = \underbrace{\bE[U(\theta_n^\eta)] - \bE[U(\theta_\infty)]}_{\mathcal{T}_1} + \underbrace{\bE[U(\theta_\infty)] - U_\star}_{\mathcal{T}_2},
\end{align*}
where $\theta_\infty \sim \pi_{\beta}$. The following proposition presents a bound for $\mathcal{T}_1$ under our assumptions.
\begin{prop}\label{prop:BoundT1} Under the assumptions of Theorem~\ref{thm:ConvRate}, we have,
\begin{align}
\bE[U(\theta_n^\eta)] - \bE[U(\theta_\infty)] &\leq \overline{C}_1^\star {\eta}^{1/2} + \overline{C}_2^\star \eta^{1/4} +\overline{C}_3^\star e^{-C_4^\star \eta n},
\end{align}
where $\overline{C}^\star_i = C^\star_i (C_m \bar{L}_1 + h_0)$ for $i = 1,2,3$ and $C^2_m = \max(C_\theta^c, C_\theta)$.
\end{prop}
\begin{proof}
See Appendix~\ref{proof:prop:BoundT1}.
\end{proof}
Next, we bound the second term $\mathcal{T}_2$ as follows. This result is fairly standard in the literature (see, e.g., \cite{raginsky2017non,gao2018global,chau2022stochastic}).
\begin{prop}\label{prop:BoundT2} \citep{raginsky2017non} Under the assumptions of Theorem~\ref{thm:ConvRate}, we have
\begin{align*}
\bE[U(\theta_\infty)] - U_\star \leq \frac{d}{2\beta} \log \left( \frac{e \overline{L}_1}{a} \left(\frac{b\beta}{d} + 1\right) \right).
\end{align*}
\end{prop}
Merging Props.~\ref{prop:BoundT1} and \ref{prop:BoundT2} leads to the bound given in Theorem~\ref{thm:NonconvexBound} which completes our proof.
\section{Applications}
In this section, we present two applications to machine learning. First, we show that the SGHMC can be used to sample from the posterior probability measure in the context of scalable Bayesian inference. We also note that our assumptions hold in a practical setting of Bayesian logistic regression. Secondly, we provide an improved generalization bound for empirical risk minimization.
\subsection{Convergence rates for scalable Bayesian inference}
Consider a prior distribution $\pi_0(\theta)$ and a likelihood $\sfp(y_i | \theta)$ for a sequence of data points $\{y_i\}_{i=1}^M$ where $M$ is the dataset size. Often, one is interested in sampling from the posterior probability distribution $p(\theta| y_{1:M}) \md\theta \propto \pi_0(\theta) \prod_{i=1}^M \sfp(y_i | \theta) \md \theta$. This is a sampling problem of the form \eqref{eq:TargetMeasure}. The SGHMC is an MCMC method to sample from the posterior measure and, therefore, our explicit convergence rates provides a guarantee for the sampling. To see this, note that the underdamped Langevin SDE with $h(\theta) = - \nabla \log p(\theta|y_{1:M})$ converges to the extended target $\overline{\pi}$ whose $\theta$-marginal of $\overline{\pi}$ is $p(\theta|y_{1:M})$, hence the underdamped Langevin SDE samples from the posterior distribution.

We note that, our setting specifically applies to cases where $M$ is too large. More precisely, note that we have
$
h(\theta) = -\nabla \log p(\theta|y_{1:M}) = -\nabla \log \pi_0(\theta) - \sum_{i=1}^M \nabla \log \sfp(y_i|\theta).
$
When $M$ is too large, evaluating $h(\theta)$ is impractical. However, one can estimate the sum in the last term in an unbiased way. To be precise, consider random indices $i_1,\ldots,i_K \sim \{1,...,M\}$ uniformly, then one can construct a stochastic gradient by using $\mathbf{u} = \{y_{i_1},\ldots,y_{i_K}\}$ and obtaining $H(\theta,\mathbf{u}) = -\log \pi_0(\theta) - \frac{M}{K} \sum_{k=1}^K \nabla \log \sfp(y_{i_k} | \theta)$. Then, we have the following corollary follows from Theorem~\ref{thm:ConvRate}.
\begin{cor}\label{corr:BayesianCorollary} Assume that the log-posterior density $\log p(\theta | y_{1:M})$, its gradient, and stochastic gradient $H(\theta,\cdot)$ satisfy the Assumptions~\ref{ass:nonnegativity}--\ref{assum:dissipativity}. Then,
\begin{align*}
W_2(\cL(\theta_n), p(\theta | y_{1:M})) &\leq C_1^\star {\eta}^{1/2} + C_2^\star \eta^{1/4} + C_3^\star e^{-C_4^\star \eta n}.
\end{align*}
where $C_1^\star, C_2^\star, C_3^\star, C_4^\star$ are finite constants.
\end{cor}

This setting becomes practical under our assumptions, e.g., for the Bayesian logistic regression example. Consider the Gaussian mixture prior $\pi_0(\theta) \propto \exp(-f_0(\theta)) = e^{-|\theta - m|^2/2} + e^{-|\theta + m|^2/2}$ and the likelihood $\sfp(\z_i | \theta) = (1/(1 + e^{-z_i^\top \theta}))^{y_i} (1 - 1/(1 + e^{-z_i^\top \theta}))^{1-y_i}$ for $\theta \in \bR^d$ and $\z_i = (z_i,y_i)$. Then, it is shown by \cite{zhang2023nonasymptotic} that the stochastic gradient $H(\theta,\mathbf{u})$ for a mini-batch in this case satisfies assumptions \ref{ass:nonnegativity}--\ref{assum:dissipativity}. In particular, our theoretical guarantee in Theorem~\ref{thm:ConvRate} and Corollary~\ref{corr:BayesianCorollary} apply to the Bayesian logistic regression case.
\subsection{A generalization bound for machine learning}
Leveraging standard results in machine learning literature, e.g., \cite{raginsky2017non}, we can prove a generalization bound for the empirical risk minimization problem. Note that, many problems in machine learning can be written as a finite-sum minimization problem as
\begin{align}\label{eq:FiniteSum}
\theta^\star \in \arg \min_{\theta\in\bR^d} U(\theta) := \frac{1}{M} \sum_{i=1}^M f(\theta,z_i).
\end{align}
Applying the result of Theorem~\ref{thm:NonconvexBound}, one can get a convergence guarantee on $\bE[U(\theta_k^\eta)] - U_\star$. However, this does not account for the so-called \textit{generalization error}. Note that, one can see the cost function in \eqref{eq:FiniteSum} as an empirical risk (expectation) minimization problem where the risk is given by $\mathbb{U}(\theta) := \int f(\theta, z) P(\md z) = \bE[f(\theta,Z)]$, where $Z \sim P(\md z)$ is an unknown probability measure where the real-world data is sampled from. Therefore, in order to bound the generalization error, one needs to bound the error $\bE[\bU(\theta_n^\eta)] - \mathbb{U}_\star$. The generalization error can be decomposed as
\begin{align*}
\bE[\bU(\theta_n^\eta)] - \mathbb{U}_\star &= \underbrace{\bE[\bU(\theta_n^\eta)] - \bE[\bU(\theta_\infty)]}_{\cB_1} + \underbrace{\bE[\bU(\theta_\infty)] - \bE[U(\theta_\infty)]}_{\cB_2} + \underbrace{\bE[U(\theta_\infty)] - \bU_\star}_{\cB_3}.
\end{align*}
In what follows, we present a series of results bounding the terms $\cB_1,\cB_2,\cB_3$. By using the results about Gibbs distributions presented in \cite{raginsky2017non}, one can bound $\cB_1$ as follows.
\begin{prop}\label{prop:GeneralizationB1} Under the assumptions of Theorem~\ref{thm:ConvRate}, we obtain
\begin{align*}
\bE[\bU(\theta_n^\eta)] - \bE[\bU(\theta_\infty)] &\leq \overline{C}_1^\star {\eta}^{1/2} + \overline{C}_2^\star \eta^{1/4} +\overline{C}_3^\star e^{-C_4^\star \eta n},
\end{align*}
where $\overline{C}^\star_i = C^\star_i (C_m \bar{L}_1 + h_0)$ for $i = 1,2,3$ and $C_m^2 = \max(C_\theta^c, C_\theta)$.
\end{prop}
The proof of Proposition~\ref{prop:GeneralizationB1} is similar to the proof of Proposition~\ref{prop:BoundT1} and indeed the rates match.

Next, we seek a bound for the term $\cB_2$. In order to prove the following result, we assume that Assumption~\ref{loclip} and Assumption~\ref{assum:dissipativity} hold uniformly in $x$, as required by, e.g., \cite{raginsky2017non}.
\begin{prop}\label{prop:GeneralizationB2} \citep{raginsky2017non} Assume that Assumptions~\ref{ass:nonnegativity},~\ref{assmp:iid} hold and Assumptions~\ref{loclip} and \ref{assum:dissipativity} hold uniformly in $x$, i.e., $|H(\theta,x)| \leq L_1' |\theta|^2 + B_1$. Then,
\begin{align*}
\bE[\bU(\theta_\infty)] - \bE[U(\theta_\infty)] \leq \frac{4 \beta c_{\textnormal{LS}}}{M} \left( \frac{L_1'}{a} (b + d/\beta) + B_1  \right),
\end{align*}
where $c_{\textnormal{LS}}$ is the constant of the logarithmic Sobolev inequality.
\end{prop}
Finally, let $\Theta^\star \in \arg \min_{\theta \in \bR} \bU(\theta)$. We note that $\cB_3$ is bounded trivially as
\begin{align}
\bE[U(\theta_\infty)] - \bU_\star &= \bE[U(\theta_\infty) - U_\star] + \bE[U_\star - U(\Theta^\star)] \leq \bE[U(\theta_\infty) - U_\star],\label{eq:GeneralizationB3}
\end{align}
which follows from the proof of Proposition~\ref{prop:BoundT2}. Finally, Proposition~\ref{prop:GeneralizationB1}, Proposition~\ref{prop:GeneralizationB2} and \eqref{eq:GeneralizationB3} leads to the following generalization bound presented as a corollary.
\begin{cor} Under the setting of Proposition~\ref{prop:GeneralizationB2}, we obtain the generalization bound for the SGHMC,
\begin{align}
\bE[\bU(\theta_n^\eta)] - \mathbb{U}_\star &\leq  \overline{C}_1^\star {\eta}^{1/2}  + \overline{C}_2^\star \eta^{1/4} +\overline{C}_3^\star e^{-C_4^\star \eta n} + \frac{4 \beta c_{\textnormal{LS}}}{M} \left( \frac{L_1'}{a} (b + d/\beta) + B_1  \right) \nonumber \\
&+ \frac{d}{2\beta} \log \left( \frac{e \overline{L}_1}{a} \left(\frac{b\beta}{d} + 1\right) \right).
\end{align}
\end{cor}
We note that this generalization bound improves that of \cite{raginsky2017non,gao2018global,chau2022stochastic} due to our improved $W_2$ bound which is reflected in Theorem~\ref{thm:NonconvexBound} and, consequently, Proposition~\ref{prop:GeneralizationB1}. In particular, while the generalization bounds of \cite{raginsky2017non} and \cite{gao2018global} grow with the number of iterations and require careful tuning between the step-size and the number of iterations, our bound decreases with the number of iterations $n$. We also note that our bound improves that of \cite{chau2022stochastic}, similar to the $W_2$ bound.
\section{Conclusions}
We have analyzed the convergence of the SGHMC recursions \eqref{eq:SGHMC1}--\eqref{eq:SGHMC2} to the extended target measure $\overline{\pi}_\beta$ in Wasserstein-2 distance which implies the convergence of the law of the iterates $\cL(\theta_n^\eta)$ to the target measure $\pi_\beta$ in $W_2$. We have proved that the error bound scales like $\mathcal{O}(\eta^{1/4})$ where $\eta$ is the step-size. This improves the existing bounds for the SGHMC significantly which are either growing with the number of iterations or include constants cannot be made to vanish by decreasing the step-size $\eta$. This bound on sampling from $\pi_\beta$ enables us to prove a stochastic global optimization result when $(\theta_n^\eta)_{n\in\bN}$ is viewed as an output of a nonconvex optimizer. We have shown that our results provide convergence rates for scalable Bayesian inference and we have particularized our results to the Bayesian logistic regression. Moreover, we have shown that our improvement of $W_2$ bounds are reflected in improved generalization bounds for the SGHMC.

\acks{\"O.~D.~A. is supported by the Lloyd’s Register Foundation Data Centric Engineering Programme and EPSRC Programme Grant EP/R034710/1 (CoSInES). S.S. acknowledges support by the Alan Turing Institute under the EPSRC grant EP/N510129/1.}

\bibliography{draft}

\newpage

\appendix
\onecolumn
\begin{huge}
\begin{center}
Appendix
\end{center}
\end{huge}
\footnotesize
\section{Constants in Theorem~\ref{thm:Eberle}}\label{app:Thm31}
The constants of Theorem~\ref{thm:Eberle} are given as follows \cite{eberle2019couplings}.
\begin{align*}
\dot{C} = 2 e^{2+\Lambda} &\frac{(1+\gamma)^2}{\min(1,\alpha)^2} \max\left(1, 4 (1 + 2\alpha + 2\alpha^2) \frac{(d+A_c)\gamma^{-1} \dot{c}^{-1}}{\min(1,R_1) \beta} \right)
\end{align*}
and
\begin{align*}
\dot{c} = \frac{\gamma}{384} \min\left(\lambda \overline{L}_1 \beta^{-1} \gamma^{-2}, \Lambda^{1/2} e^{-\Lambda} \overline{L}_1 \beta^{-1} \gamma^{-2}, \Lambda^{1/2} e^{-\Lambda}\right)
\end{align*}
where
\begin{align*}
\Lambda &= \overline{L}_1 R_1^2 / 8 = \frac{12}{5} (1 + 2\alpha + 2\alpha^2)(d+A_c)\overline{L}_1 \beta^{-1} \gamma^{-2} \lambda ^{-1} (1 - 2\lambda)^{-1},
\end{align*}
and $\alpha \in (0,\infty)$.
\section{Additional lemmata}
We prove the following lemma adapted from \cite{raginsky2017non}.
\begin{lem}\label{lem:UpperLowerQuadBoundForU} For all $\theta \in \bR^d$,
\begin{align*}
\frac{a}{3} |\theta|^2 - \frac{b}{2} \log 3 \leq U(\theta) \leq u_0 + \frac{\bar{L}_1}{2} |\theta|^2 + {h_0 |\theta|}.
\end{align*}
where $u_0 = U(0)$ and $\bar{L}_1 = L_1 \bE[(1+|X_0|)^\rho]$.
\end{lem}
\begin{proof}
See Appendix~\ref{proof:lem:UpperLowerQuadBoundForU}.
\end{proof}

\begin{lem}\label{lem:magicLemma} Let $\mathcal{G},\mathcal{H} \subset \cF$ be sigma-algebras. Consider two $\bR^d$-valued random vectors, denoted $X,Y$, in $L^p$ with $p \geq 1$ such that $Y$ is measurable w.r.t. $\mathcal{H} \vee \mathcal{G}$. Then,
\begin{align*}
\bE^{1/p}[ | X - \bE[ X | \mathcal{H} \vee \mathcal{G}] |^p | \mathcal{G}] \leq 2 \bE^{1/p}[ | X - Y |^p | \mathcal{G}].
\end{align*}
\end{lem}
\begin{proof}
See Lemma~6.1 in \cite{chau2019fixed}.
\end{proof}

\begin{lem}\label{lem:boundedVariance} Let Assumption \ref{ass:nonnegativity}, \ref{assmp:iid}, \ref{loclip} and \ref{assum:dissipativity} hold. For any $k = 1, \dots, K+1$ where $K+1 \leq T$, we obtain
\begin{align*}
\sup_{n\in\bN} \sup_{t  \in [nT, (n+1)T]}\E\left[\left|h(\bar{\zeta}_t^{\eta,n}) - H(\bar{\zeta}_t^{\eta,n},X_{nT+k})\right|^2\right] \leq \sigma_H
\end{align*}
where
\begin{align*}
\sigma_H := 8 L_2^2 \sigma_Z (1 + C_\zeta) < \infty,
\end{align*}
where $\sigma_Z = \bE[(1 + |X_0| + | \bE[X_0]|)^{2\rho} |X_0 - \bE[X_0]|^2] < \infty$.
\end{lem}
\begin{proof}
See Appendix~\ref{proof:lem:boundedVariance}.
\end{proof}

\begin{lem}\label{lem:distanceToAuxiliary} Under Assumptions \ref{ass:nonnegativity}--\ref{assum:dissipativity}
\begin{align*}
\bE\left[ \left| \overline{V}_{\floor{t}}^\eta - \overline{V}_t^\eta \right|^2\right] \leq {\sigma}_V \eta,
\end{align*}
where
\begin{align*}
{\sigma}_V = 4 \eta \gamma^2 C_v + 4 \eta (\widetilde{L}_1 C_\theta + \widetilde{C}_1) + 4 \gamma \beta^{-1} d.
\end{align*}
Moreover,
\begin{align*}
\bE\left[\left| \overline{\theta}_\floor{t}^\eta - \overline{\theta}_t^\eta \right|^2\right] \leq \eta^2 C_v.
\end{align*}
\end{lem}
\begin{proof}
See Appendix~\ref{proof:lem:distanceToAuxiliary}.
\end{proof}
Next, it is shown that a key assumption appearing in \cite{eberle2019couplings} holds.
\begin{lem}\label{lem:LowerBoundOfInnerProd} There exist constants $A_c \in (0,\infty)$ and $\lambda \in (0,1/4]$ such that
\begin{align}\label{eq:DissipativityLowerBound}
\langle \theta, h(\theta)\rangle \geq 2 \lambda \left( U(\theta) + \gamma^2 |\theta |^2 / 4 \right) - 2 A_c/\beta
\end{align}
for all $x \in \bR^d$.
\end{lem}
\begin{proof}
See Appendix~\ref{proof:lem:LowerBoundOfInnerProd}.
\end{proof}

\section{Proofs of the preliminary results}\label{app:sec:proof:prel_res}

\subsection{Proof of Lemma~\ref{lem:DiscreteTimeLemma}}\label{proof:lem:DiscreteTimeLemma}
For this proof, we use the Lyapunov function defined by \cite{eberle2019couplings} and follow a similar proof presented in \cite{gao2018global}. We first define the Lyapunov function as
\begin{align*}
\mathcal{V}(\theta,v) = \beta U(\theta) + \frac{\beta}{4} \gamma^2 \left( \|\theta + \gamma^{-1} v\|^2 + \|\gamma^{-1} v\|^2 - \lambda \|\theta\|^2 \right).
\end{align*}
Next, we will use this Lyapunov function to show that the second moments of the processes $(V^\eta_n)_{n\in\bN}$ and $(\theta^\eta_n)_{n\in\bN}$ are finite.

We start by defining
\begin{align}\label{eq:DefnL2funct}
M_2(k) = \bE \mathcal{V}(\theta_k^\eta, V_k^\eta)/\beta = \bE \left[ U(\theta_k^\eta) + \frac{\gamma^2}{4} \left( |\theta_k^\eta + \gamma^{-1} V_k^\eta|^2 + |\gamma^{-1} V_k^\eta|^2 - \lambda |\theta_k^\eta|^2 \right) \right].
\end{align}
Recall our discrete-time recursions \eqref{eq:SGHMC1}--\eqref{eq:SGHMC2}
\begin{align*}
{V}_{k+1}^\eta &= (1 - \eta\gamma) {V}_k^\eta - \eta H({\theta}_k^\eta, X_{k+1}) + \sqrt{2\eta\gamma\beta^{-1}} \xi_{k+1}, \\
\theta_{k+1}^\eta &= \theta_k^\eta + \eta V_k^\eta, \quad \quad \quad \quad \theta_0^\eta = \theta_0, \quad \quad V_0^\eta = v_0,
\end{align*}
where $(\xi_k)_{k\in\bN}$ is a sequence of i.i.d. standard Normal random variables. Consequently, we have the equality
\begin{align*}
\bE\left[| {V}_{k+1}^\eta|^2 \right] &= \bE\left[\left| (1-\gamma\eta) V_k^\eta - \eta H(\theta_k^\eta,X_{k+1}) \right|^2 \right] + 2\gamma\eta \beta^{-1} d, \\
&= (1 - \gamma\eta)^2 \bE\left[ \left| V_k^\eta \right|^2\right] - 2 \eta (1 - \gamma\eta) \bE \left[ \langle V_k^\eta, h(\theta_k^\eta)\rangle\right] + \eta^2 \bE\left[ \left| H(\theta_k^\eta, X_{k+1})\right|^2\right] + 2\gamma\eta \beta^{-1} d,
\end{align*}
which immediately leads to
\begin{align}\label{eq:BoundVkSquare}
\bE\left[| {V}_{k+1}^\eta|^2 \right] &\leq (1 - \gamma\eta)^2 \bE\left[ \left| V_k^\eta \right|^2\right] - 2 \eta (1 - \gamma\eta) \bE\left[ \langle V_k^\eta, h(\theta_k^\eta)\rangle\right] + \eta^2 \left[ \widetilde{L}_1 \bE\left[ \left| \theta_k^\eta\right|^2\right] + \widetilde{C}_1 \right] + 2\gamma\eta\beta^{-1} d,
\end{align}
where
\begin{align}\label{eq:LtildeCtilde}
\widetilde{L}_1 = 2 L_1^2 C_\rho \quad \quad \textnormal{and} \quad \quad \widetilde{C}_1 = 4 L_2^2 C_\rho + 4 H_0^2.
\end{align}
Next, we note that
\begin{align}\label{eq:EqualityThetakSquare}
\bE \left| \theta_{k+1}^\eta \right|^2 = \bE | \theta_k^\eta|^2 + 2 \eta \bE \langle \theta_k^\eta, V_k^\eta \rangle + \eta^2 \bE \left| V_k^\eta \right|^2.
\end{align}
Recall $h := \nabla U$ and note also that
\begin{align*}
U(\theta_{k+1}^\eta) = U(\theta_k^\eta + \eta V_k^\eta) = U(\theta_k^\eta) + \int_0^1 \langle h(\theta_k^\eta + \tau\eta V_k^\eta), \eta V_k^\eta\rangle \md \tau,
\end{align*}
which suggests
\begin{align*}
\left| U(\theta_{k+1}^\eta) - U(\theta_k^\eta) - \langle h(\theta_k^\eta), \eta V_k^\eta) \rangle \right| &= \left| \int_0^1 \langle h(\theta_k^\eta + \tau \eta V_k^\eta) - {h(\theta_k^\eta)}, \eta V_k^\eta \rangle \md \tau \right|, \\
&\leq \int_0^1 \left| h(\theta_k^\eta + \tau \eta V_k^\eta) - h(\theta_k^\eta) \right| \left| \eta V_k^\eta\right| \md \tau, \\
&\leq \frac{1}{2} L_1 C_\rho \eta^2 \left| V_k^\eta \right|^2,
\end{align*}
where the second line follows from the Cauchy-Schwarz inequality and the final line follows from \eqref{eq:LipschitzGradient}. Finally we obtain
\begin{align}\label{eq:BoundExpFuncValueDecrease}
\bE U(\theta_{k+1}^\eta) - \bE U(\theta_k^\eta) \leq \eta \bE \langle h(\theta_k^\eta), V_k^\eta \rangle + \frac{1}{2} L_1 C_\rho \eta^2 \bE \left|V_k^\eta \right|^2.
\end{align}
Next, we continue computing
\begin{align}
\bE \left| \theta_{k+1}^\eta + \gamma^{-1} V_{k+1}^\eta \right|^2 &= \bE \left| \theta_k^\eta + \gamma^{-1} V_k^\eta - \eta \gamma^{-1} H(\theta_k^\eta, X_{k+1})\right|^2 + 2 \gamma^{-1} \beta^{-1} \eta d, \nonumber\\
&= \bE \left| \theta_k^\eta + \gamma^{-1} V_k^\eta \right|^2 - 2 \eta \gamma^{-1} \bE \langle \theta_k^\eta + \gamma^{-1} V_k^\eta, h(\theta_k^\eta)\rangle \nonumber \\
& + \eta^2 \gamma^{-2} \bE \left| H(\theta_k^\eta, X_{k+1})\right|^2 + 2 \gamma^{-1} \eta \beta^{-1} d, \nonumber \\
&\leq \bE \left| \theta_k^\eta + \gamma^{-1} V_k^\eta \right|^2 - 2 \eta \gamma^{-1} \bE \langle \theta_k^\eta + \gamma^{-1} V_k^\eta, h(\theta_k^\eta)\rangle \nonumber\\
& + \eta^2 \gamma^{-2} (\widetilde{L}_1 \bE\left| \theta_k^\eta\right|^2 + \widetilde{C}_1) + 2 \gamma^{-1} \eta \beta^{-1} d. \label{eq:LyapunovTermDecrease}
\end{align}
where $\widetilde{L}_1$ and $\widetilde{C}_1$ is defined as in \eqref{eq:LtildeCtilde}. Next, combining \eqref{eq:BoundVkSquare}, \eqref{eq:EqualityThetakSquare}, \eqref{eq:BoundExpFuncValueDecrease}, \eqref{eq:LyapunovTermDecrease},
\begin{align*}
M_2(k&+1) - M_2(k) \\ &= \bE\left[ U(\theta_{k+1}^\eta) - U(\theta_k^\eta)\right] + \frac{\gamma^2}{4} \left( \bE \left| \theta_{k+1}^\eta + \gamma^{-1} V_{k+1}^\eta \right|^2 - \bE \left| \theta_k^\eta + \gamma^{-1} V_k^\eta\right|^2 \right) \\
&+ \frac{1}{4} \left( \bE\left| V_{k+1}^\eta\right|^2 - \bE\left|V_k^\eta\right|^2\right) - \frac{\gamma^2 \lambda}{4} \left( \bE\left|\theta_{k+1}^\eta\right|^2 - \bE\left|\theta_k^\eta\right|^2 \right), \\
&\leq \eta \bE \langle h(\theta_k^\eta),V_k^\eta\rangle + \frac{L_1 C_\rho \eta^2}{2} \bE\left| V_k^\eta \right|^2 \\
&+ \frac{\gamma^2}{4} \left(-2\eta\gamma^{-1} \bE \langle \theta_k^\eta + \gamma^{-1} V_k^\eta, h(\theta_k^\eta)\rangle + \eta^2\gamma^{-2} \left(\widetilde{L}_1 \bE\left|\theta_k^\eta\right|^2 + \widetilde{C}_1\right) +  2 \gamma^{-1} \beta^{-1} \eta d \right) \\
&+ \frac{1}{4} \left( (-2\gamma\eta + \gamma^2\eta^2) \bE\left|V_k^\eta\right|^2 - 2 \eta (1 - \gamma \eta) \bE \langle V_k^\eta, h(\theta_k^\eta)\rangle + \eta^2 \left(\widetilde{L}_1 \bE|\theta_k^\eta|^2 + \widetilde{C}_1\right) + 2\gamma \eta \beta^{-1} d \right) \\
&- \frac{\gamma^2 \lambda}{4} \left( 2 \eta \bE \langle \theta_k^\eta, V_k^\eta\rangle + \eta^2 \bE |V_k^\eta|^2 \right),\\
&=  - \frac{\eta\gamma}{2} \bE\langle \theta_k^\eta, h(\theta_k^\eta)\rangle + \frac{\gamma\eta^2}{2} \bE \langle h(\theta_k^\eta),V_k^\eta\rangle + \left(\frac{{L}_1 C_\rho \eta^2}{2} + \frac{\eta^2 \gamma^2}{4} - \frac{\gamma\eta}{2} - \frac{\gamma^2\eta^2 \lambda}{4}\right) \bE\left| V_k^\eta \right|^2  \\
&+ \frac{\eta^2 \widetilde{L}_1}{2} \bE \left| \theta_k^\eta\right|^2 - \frac{\gamma^2\eta\lambda}{2} \bE\langle \theta_k^\eta,V_k^\eta\rangle + \frac{\widetilde{C_1}\eta^2}{2} + {\gamma\eta\beta^{-1}d}, \\
&\leq - \eta\gamma\lambda \bE U(\theta_k^\eta) - \frac{\lambda \gamma^3 \eta}{4} \bE |\theta_k^\eta|^2 + A_c \eta\gamma \beta^{-1} + \frac{\gamma\eta^2}{2} \bE \langle h(\theta_k^\eta),V_k^\eta\rangle + \frac{\eta^2 \widetilde{L}_1}{2} \bE \left| \theta_k^\eta\right|^2 \\& + \left(\frac{{L}_1 C_\rho \eta^2}{2} + \frac{\eta^2 \gamma^2}{4} - \frac{\gamma\eta}{2} - \frac{\gamma^2\eta^2 \lambda}{4}\right) \bE\left| V_k^\eta \right|^2  - \frac{\gamma^2\eta\lambda}{2} \bE\langle \theta_k^\eta,V_k^\eta\rangle + \frac{\widetilde{C_1}\eta^2}{2} + {\gamma\eta\beta^{-1}d}.
\end{align*}
where the last line is obtained using \eqref{eq:DissipativityLowerBound}. Next, using the fact that $0 < \lambda \leq 1/4$ and the form of the Lyapunov function \eqref{eq:DefnL2funct}, we obtain
\begin{align*}
- \frac{\gamma}{2} \bE \langle \theta_k^\eta,V_k^\eta\rangle \leq - M_2(k) + \bE U(\theta_k^\eta) + \frac{\gamma^2}{4} \bE |\theta_k^\eta|^2 + \frac{1}{2} \bE |V_k^\eta|^2.
\end{align*}
Using this, we can obtain
\begin{align*}
&M_2(k+1) - M_2(k) \leq {A_c \eta\gamma}{\beta^{-1}} + \frac{\gamma\eta^2}{2} \bE \langle h(\theta_k^\eta),V_k^\eta\rangle + \frac{\eta^2 \widetilde{L}_1}{2} \bE \left| \theta_k^\eta\right|^2 + {\gamma\eta\beta^{-1}d} \\& + \left(\frac{{L}_1 C_\rho \eta^2}{2} + \frac{\eta^2 \gamma^2}{4} - \frac{\gamma\eta}{2} - \frac{\gamma^2\eta^2 \lambda}{4} + \frac{\gamma\lambda\eta}{2} \right) \bE\left| V_k^\eta \right|^2  + \frac{\widetilde{C_1}\eta^2}{2} - \gamma\lambda\eta M_2(k).
\end{align*}
Next, reorganizing and using $\langle a,b\rangle \leq (|a|^2 + |b|^2)/2$
\begin{align*}
M_2(k+1) &\leq (1 - \gamma\lambda\eta) M_2(k) + {A_c \eta\gamma}{\beta^{-1}} + \frac{\eta^2 \widetilde{L}_1}{2} \bE \left| \theta_k^\eta\right|^2 + \frac{\widetilde{C_1}\eta^2}{2} + {\gamma\eta\beta^{-1}d} \\&
+ \left(\frac{{L}_1 C_\rho \eta^2}{2} + \frac{\eta^2 \gamma^2}{4} - \frac{\gamma\eta}{2} - \frac{\gamma^2\eta^2 \lambda}{4} + \frac{\gamma\lambda\eta}{4} + \frac{\gamma\eta^2}{4} \right) \bE\left| V_k^\eta \right|^2  + \frac{\gamma\eta^2}{4} \bE |h(\theta_k^\eta)|^2, \\
&\leq (1 - \gamma\lambda\eta) M_2(k) + {A_c \eta\gamma}{\beta^{-1}} + \eta^2 \left(\frac{\widetilde{L}_1}{2} + \frac{\gamma}{2} {L^2_1} C_\rho^2 \right) \bE \left| \theta_k^\eta\right|^2 + \frac{\widetilde{C_1}\eta^2}{2} + {\gamma\eta\beta^{-1}d} \\& + \eta^2 \left(\frac{{L}_1 C_\rho}{2} + \frac{ \gamma^2}{4} - \frac{\gamma^2 \lambda}{4} + \frac{\gamma}{4} \right) \bE\left| V_k^\eta \right|^2 + \frac{\gamma \eta^2 h_0^2}{2},
\end{align*}
where the last inequality follows since $\lambda \leq 1/4$ and \eqref{eq:BoundedGradh}. We note that
\begin{align*}
\mathcal{V}(\theta,v) \geq \max\left\lbrace \frac{1}{8} (1 - 2\lambda) \beta \gamma^2 |\theta|^2, \frac{\beta}{4} (1 - 2\lambda) |v|^2 \right\rbrace,
\end{align*}
which implies by the definition of $M_2(k)$ that
\begin{align}
M_2(k) &\geq \max\left\lbrace \frac{1}{8} (1 - 2\lambda) \gamma^2 \bE|\theta_k^\eta|^2, \frac{1}{4} (1 - 2\lambda) \bE |V_k^\eta|^2 \right\rbrace, \nonumber \\
&\geq \frac{1}{16} (1 - 2\lambda) \gamma^2 \bE|\theta_k^\eta|^2 + \frac{1}{8} (1 - 2\lambda) \bE |V_k^\eta|^2, \label{eq:UpperBoundOfBothMoments}
\end{align}
since $\max\{x,y\} \geq (x+y) / 2$ for any $x,y > 0$. Therefore, we obtain
\begin{align*}
M_2(k+1) &\leq (1 - \gamma\lambda\eta + K_1 \eta^2) M_2(k) + K_2 \eta^2 + K_3 \eta
\end{align*}
where
\begin{align*}
K_1 := \max\left\lbrace \frac{\frac{{L}_1 C_\rho}{2} + \frac{ \gamma^2}{4} - \frac{\gamma^2 \lambda}{4} + \frac{\gamma}{4}}{\frac{1}{8} (1 - 2\lambda)},  \frac{\frac{\widetilde{L}_1}{2} + \frac{\gamma}{2} {L^2_1} C_\rho^2}{\frac{1}{16} (1 - 2\lambda) \gamma^2} \right\rbrace
\end{align*}
and
\begin{align*}
K_2 = \frac{\widetilde{C_1} + \gamma h_0^2}{2} \quad \textnormal{and} \quad K_3 = (A_c + d) \gamma \beta^{-1}.
\end{align*}
For $0 < \eta \leq \min\left\lbrace \frac{K_3}{K_2}, \frac{\gamma\lambda}{2 K_1}, \frac{2}{\gamma\lambda}\right\rbrace$, we obtain
\begin{align*}
M_2(k+1) \leq \left(1 - \frac{{\gamma\lambda\eta}}{2}\right) M_2(k) + 2 K_3\eta
\end{align*}
which implies
\begin{align*}
M_2(k) \leq M_2(0) + \frac{4}{\gamma\lambda} K_3.
\end{align*}
Combining this with \eqref{eq:UpperBoundOfBothMoments} gives the result.

\subsection{Proof of Lemma~\ref{lem:zetaprocessmoment}}\label{proof:lem:zetaprocessmoment}
We recall that $\overline{\zeta}_t^{\eta,n}$ is the Langevin diffusion started at $\theta^\eta_{nT}$ and run until $t \in (nT,(n+1)T)$. First notice that Lemma~\ref{lem:ContTimeLemma} implies
\begin{align}
\sup_{t \in (nT, (n+1) T)} \bE | \zeta_t^{\eta,n}|^2 \leq \frac{\bE \mathcal{V}(\theta_{nT}^\eta, V_{nT}^\eta) + \frac{4 (d+A_c)}{\lambda}}{\frac{1}{8} (1 - 2\lambda) \beta \gamma^2},
\end{align}
which, noting that, $\bE \mathcal{V}(\theta_{nT}^\eta, V_{nT}^\eta) = \beta M_2(nT)$, implies
\begin{align*}
\sup_{t \in (nT, (n+1) T)} \bE | \zeta_t^{\eta,n}|^2 \leq \frac{\beta M_2(0) + \frac{8(d+A_c)}{\lambda}}{\frac{1}{8} (1 - 2\lambda) \beta \gamma^2}.
\end{align*}
Substituting $M_2(0)$ gives
\begin{align*}
\sup_{t \in (nT, (n+1) T)} \bE | \zeta_t^{\eta,n}|^2 \leq C_\zeta := \frac{\int_{\bR^{2d}} \mathcal{V}(\theta,v) \mu_0(\md \theta,\md v) + \frac{8(d + A_c)}{\lambda}}{\frac{1}{8} (1-2\lambda)\beta\gamma^2}.
\end{align*}


\subsection{Proof of Lemma~\ref{lem:SquareContraction}}\label{proof:lem:SquareContraction}
We need to obtain the contraction of $\bE\cV^2(\theta_k^\eta, V_k^\eta)$. Recall again the Lyapunov function defined by \cite{eberle2019couplings}
\begin{align*}
\mathcal{V}(\theta,v) = \beta U(\theta) + \frac{\beta}{4} \gamma^2 \left( \|\theta + \gamma^{-1} v\|^2 + \|\gamma^{-1} v\|^2 - \lambda \|\theta\|^2 \right).
\end{align*}
In this proof, we follow a similar strategy to the proof of Lemma~\ref{lem:DiscreteTimeLemma} and approaches of \cite{gao2018global} and \cite{chau2022stochastic}.
We define the notation $\cV_k = \mathcal{V}(\theta_k^\eta, V_k^\eta)$ and note
\begin{align}\label{eq:DefnMfunct}
\cV_k/\beta := \mathcal{V}(\theta_k^\eta, V_k^\eta)/\beta = U(\theta_k^\eta) + \frac{\gamma^2}{4} \left( |\theta_k^\eta + \gamma^{-1} V_k^\eta|^2 + |\gamma^{-1} V_k^\eta|^2 - \lambda |\theta_k^\eta|^2 \right).
\end{align}
Recall our discrete-time recursions \eqref{eq:SGHMC1}--\eqref{eq:SGHMC2}
\begin{align*}
{V}_{k+1}^\eta &= (1 - \eta\gamma) {V}_k^\eta - \eta H({\theta}_k^\eta, X_{k+1}) + \sqrt{2\eta\gamma\beta^{-1}} \xi_{k+1}, \\
\theta_{k+1}^\eta &= \theta_k^\eta + \eta V_k^\eta, \quad \quad \quad \quad \theta_0^\eta = \theta_0, \quad \quad V_0^\eta = v_0,
\end{align*}
where $(\xi_k)_{k\in\bN}$ is a sequence of i.i.d. Normal random variables. We first define $\Delta_k^1=(1-\eta\gamma)V_k^\eta~-\eta~H(\theta_k^\eta, X_{k+1})$ and write
\begin{align}
| {V}_{k+1}^\eta|^2 &= |\Delta_k^1|^2 + 2 \sqrt{2 \eta \gamma \beta^{-1}} \langle \Delta_k^1, \xi_{k+1} \rangle + 2\gamma\eta \beta^{-1} |\xi_{k+1}|^2, \nonumber \\
&= \left| (1-\gamma\eta) V_k^\eta - \eta H(\theta_k^\eta,X_{k+1}) \right|^2 + 2 \sqrt{2 \eta \gamma \beta^{-1}} \langle \Delta_k^1, \xi_{k+1} \rangle + 2\gamma\eta \beta^{-1} |\xi_{k+1}|^2, \nonumber \\
&= (1 - \gamma\eta)^2 \left| V_k^\eta \right|^2 - 2 \eta (1 - \gamma\eta)  \langle V_k^\eta, H(\theta_k^\eta,X_{k+1})\rangle \nonumber \\ &+ \eta^2 \left| H(\theta_k^\eta, X_{k+1})\right|^2 + 2 \sqrt{2 \eta \gamma \beta^{-1}} \langle \Delta_k^1, \xi_{k+1} \rangle + 2\gamma\eta \beta^{-1} |\xi_{k+1}|^2. \label{eq:AlmostSureEqualityVkSquare}
\end{align}
Next, we note that
\begin{align}\label{eq:AlmostSureEqualityThetakSquare}
\left| \theta_{k+1}^\eta \right|^2 = | \theta_k^\eta|^2 + 2 \eta \langle \theta_k^\eta, V_k^\eta \rangle + \eta^2 \left| V_k^\eta \right|^2.
\end{align}
Recall $h := \nabla U$ and note also that
\begin{align*}
U(\theta_{k+1}^\eta) = U(\theta_k^\eta + \eta V_k^\eta) = U(\theta_k^\eta) + \int_0^1 \langle h(\theta_k^\eta + \tau\eta V_k^\eta), \eta V_k^\eta\rangle \md \tau,
\end{align*}
which suggests
\begin{align*}
\left| U(\theta_{k+1}^\eta) - U(\theta_k^\eta) - \langle h(\theta_k^\eta), \eta V_k^\eta) \rangle \right| &= \left| \int_0^1 \langle h(\theta_k^\eta + \tau \eta V_k^\eta) - {h(\theta_k^\eta)}, \eta V_k^\eta \rangle \md \tau \right|, \\
&\leq \int_0^1 \left| h(\theta_k^\eta + \tau \eta V_k^\eta) - h(\theta_k^\eta) \right| \left| \eta V_k^\eta\right| \md \tau, \\
&\leq \frac{1}{2} L_1 C_\rho \eta^2 \left| V_k^\eta \right|^2,
\end{align*}
where, similarly as in the previous proof, the second line follows from the Cauchy-Schwarz inequality and the final line follows from \eqref{eq:LipschitzGradient}. Finally we arrive at
\begin{align}\label{eq:AlmostSureBoundExpFuncValueDecrease}
U(\theta_{k+1}^\eta) - U(\theta_k^\eta) \leq \eta  \langle h(\theta_k^\eta), V_k^\eta \rangle + \frac{1}{2} L_1 C_\rho \eta^2 \left|V_k^\eta \right|^2.
\end{align}
Next, we note that $\Delta_k^2 = \theta_k^\eta + \gamma^{-1} V_k^\eta - \eta \gamma^{-1} H(\theta_k^\eta, X_{k+1})$ and continue computing
\begin{align}
\left| \theta_{k+1}^\eta + \gamma^{-1} V_{k+1}^\eta \right|^2 &= \left| \theta_k^\eta + \gamma^{-1} V_k^\eta - \eta \gamma^{-1} H(\theta_k^\eta, X_{k+1})\right|^2 + 2\sqrt{2\gamma^{-1}\beta^{-1} \eta} \langle \Delta_k^2, \xi_{k+1}\rangle + 2 \gamma^{-1} \beta^{-1} \eta |\xi_{k+1}|^2, \nonumber\\
&= \left| \theta_k^\eta + \gamma^{-1} V_k^\eta \right|^2 - 2 \eta \gamma^{-1} \langle \theta_k^\eta + \gamma^{-1} V_k^\eta, H(\theta_k^\eta,X_{k+1})\rangle \nonumber \\
& + 2\sqrt{2\gamma^{-1}\beta^{-1} \eta} \langle \Delta_k^2, \xi_{k+1}\rangle  + \eta^2 \gamma^{-2}  \left| H(\theta_k^\eta, X_{k+1})\right|^2 + 2 \gamma^{-1} \eta \beta^{-1} |\xi_{k+1}|^2, \label{eq:AlmostSureLyapunovTerm}
\end{align}
Next, combining \eqref{eq:AlmostSureEqualityVkSquare}, \eqref{eq:AlmostSureEqualityThetakSquare}, \eqref{eq:AlmostSureBoundExpFuncValueDecrease}, \eqref{eq:AlmostSureLyapunovTerm},
\begin{align*}
&\frac{\cV_{k+1} - \cV_k}{\beta} = \left( U(\theta_{k+1}^\eta) - U(\theta_k^\eta)\right) + \frac{\gamma^2}{4} \left( \left| \theta_{k+1}^\eta + \gamma^{-1} V_{k+1}^\eta \right|^2 - \left| \theta_k^\eta + \gamma^{-1} V_k^\eta\right|^2 \right) \\
&+ \frac{1}{4} \left( \left| V_{k+1}^\eta\right|^2 - \left|V_k^\eta\right|^2\right) - \frac{\gamma^2 \lambda}{4} \left( \left|\theta_{k+1}^\eta\right|^2 - \left|\theta_k^\eta\right|^2 \right), \\
&\leq  \eta  \langle h(\theta_k^\eta), V_k^\eta \rangle + \frac{1}{2} L_1 C_\rho \eta^2 \left|V_k^\eta \right|^2 + \frac{\gamma^2}{4} \left( - 2 \eta \gamma^{-1} \langle \theta_k^\eta + \gamma^{-1} V_k^\eta, H(\theta_k^\eta,X_{k+1})\rangle \right. \\ & \left.+ \eta^2 \gamma^{-2}  \left| H(\theta_k^\eta, X_{k+1})\right|^2 + 2 \gamma^{-1} \eta \beta^{-1} |\xi_{k+1}|^2\right) + \frac{1}{4} \left( (-2\gamma\eta + \gamma^2\eta^2) \left| V_k^\eta \right|^2 - 2 \eta (1 - \gamma\eta)  \langle V_k^\eta, H(\theta_k^\eta,X_{k+1})\rangle \right. \\ & \left.+ \eta^2 \left| H(\theta_k^\eta, X_{k+1})\right|^2 + 2\gamma\eta \beta^{-1} |\xi_{k+1}|^2\right) - \frac{\gamma^2\lambda}{4} \left(2 \eta \langle \theta_k^\eta, V_k^\eta \rangle + \eta^2 \left| V_k^\eta \right|^2 \right) + \Sigma_k, \\
&\leq - \frac{\eta\gamma}{2} \langle \theta_k^\eta, H(\theta_k^\eta, X_{k+1})\rangle + \eta \langle h(\theta_k^\eta), V_k^\eta \rangle + \left( \frac{L_1 C_\rho \eta^2}{2} - \frac{\gamma^2 \lambda \eta^2}{4} - \frac{\gamma\eta}{2} + \frac{\gamma^2\eta^2}{4} \right) |V_k^\eta|^2 + \frac{\eta^2}{2} | H(\theta_k^\eta, X_{k+1})|^2 \\
&+ \eta \gamma \beta^{-1} | \xi_{k+1} |^2 + \left( \frac{\gamma \eta^2}{2} -  \eta \right) \langle V_k^\eta, H(\theta_k^\eta, X_{k+1})\rangle - \frac{\gamma^2 \lambda \eta}{2} \langle \theta_k^\eta, V_k^\eta \rangle + \Sigma_k
\end{align*}
where
\begin{align*}
\Sigma_k = \frac{\gamma^2}{2} \sqrt{2\gamma^{-1}\beta^{-1} \eta} \langle \Delta_k^2, \xi_{k+1}\rangle + \frac{1}{2} \sqrt{2 \eta \gamma \beta^{-1}} \langle \Delta_k^1, \xi_{k+1}\rangle.
\end{align*}
Next, using the fact that $0 < \lambda \leq 1/4$ and the form of the Lyapunov function \eqref{eq:DefnMfunct}, we obtain
\begin{align*}
- \frac{\gamma}{2} \langle \theta_k^\eta,V_k^\eta\rangle \leq - \frac{\cV_k}{\beta} + U(\theta_k^\eta) + \frac{\gamma^2}{4}  |\theta_k^\eta|^2 + \frac{1}{2}  |V_k^\eta|^2.
\end{align*}
Using this and merging some terms, we obtain
\begin{align}
\frac{\cV_{k+1} - \cV_k}{\beta} &\leq - \frac{\eta\gamma}{2} \langle \theta_k^\eta, H(\theta_k^\eta, X_{k+1})\rangle + \eta \langle h(\theta_k^\eta) - H(\theta_k^\eta,X_{k+1}), V_k^\eta \rangle \nonumber \\ &+ \left( \frac{L_1 C_\rho \eta^2}{2} - \frac{\gamma^2 \lambda \eta^2}{4} - \frac{\gamma\eta}{2} + \frac{\gamma^2\eta^2}{4} \right) |V_k^\eta|^2 + \frac{\eta^2}{2} | H(\theta_k^\eta, X_{k+1})|^2 + \eta \gamma \beta^{-1} | \xi_{k+1} |^2 \nonumber \\&+ \frac{\gamma \eta^2}{2} \langle V_k^\eta, H(\theta_k^\eta, X_{k+1})\rangle - \gamma\lambda\eta \frac{\cV_k}{\beta} + \gamma\lambda\eta U(\theta_k^\eta) + \frac{\gamma^3 \lambda\eta}{4} |\theta_k^\eta|^2 + \frac{\gamma\lambda\eta}{2} | V_k^\eta |^2 + \Sigma_k. \label{eq:eqref_diff_divided_by_beta}
\end{align}
Next, by reorganizing {and using $\langle a, b\rangle \leq (|a|^2 + |b|^2)/2$}, we arrive at
\begin{align}
&\frac{\cV_{k+1}}{\beta} \leq (1 - \gamma\lambda\eta) \frac{\cV_k}{\beta} - \frac{\eta\gamma}{2} \langle \theta_k^\eta, h(\theta_k^\eta)\rangle + \frac{\eta\gamma}{2} \langle \theta_k^\eta, h(\theta_k^\eta) - H(\theta_k^\eta, X_{k+1})\rangle + \eta \langle h(\theta_k^\eta) - H(\theta_k^\eta,X_{k+1}), V_k^\eta \rangle \nonumber
\\& + \left( \frac{L_1 C_\rho \eta^2}{2} - \frac{\gamma^2 \lambda \eta^2}{4} - \frac{\gamma\eta}{2} + \frac{\gamma^2\eta^2}{4} + \frac{\gamma \lambda \eta}{2} \right) |V_k^\eta|^2  + \frac{\eta^2}{2} | H(\theta_k^\eta, X_{k+1})|^2 + \eta \gamma \beta^{-1} | \xi_{k+1} |^2 + \gamma\lambda\eta U(\theta_k^\eta) + \frac{\gamma^3 \lambda\eta}{4} |\theta_k^\eta|^2 \nonumber \\&+ \Sigma_k. \label{eqref:v_k_div_by_beta}
\end{align}
Using \eqref{eq:DissipativityLowerBound},  and $\lambda \leq 1/4$, we obtain
\begin{align*}
&\frac{\cV_{k+1}}{\beta} \leq (1 - \gamma\lambda\eta) \frac{\cV_k}{\beta} - \eta\gamma\lambda U(\theta_k^\eta) - \frac{\gamma^3 \eta \lambda}{4} |\theta_k^\eta|^2  + \eta \gamma A_c \beta^{-1} + \frac{\eta\gamma}{2} \langle \theta_k^\eta, h(\theta_k^\eta) - H(\theta_k^\eta, X_{k+1})\rangle
\\& + \eta \langle h(\theta_k^\eta) - H(\theta_k^\eta,X_{k+1}), V_k^\eta \rangle + \left( \frac{L_1 C_\rho \eta^2}{2} - \frac{\gamma^2 \lambda \eta^2}{4} - \frac{\gamma\eta}{2} + \frac{\gamma^2\eta^2}{4} + \frac{\gamma \lambda \eta}{2} \right) |V_k^\eta|^2  + \frac{\eta^2}{2} | H(\theta_k^\eta, X_{k+1})|^2 \\ &+ \eta \gamma \beta^{-1} | \xi_{k+1} |^2
 + \gamma\lambda\eta U(\theta_k^\eta) + \frac{\gamma^3 \lambda\eta}{4}  |\theta_k^\eta|^2 + \Sigma_k,
\\
&= (1 - \gamma\lambda\eta) \frac{\cV_k}{\beta} + \eta \gamma A_c \beta^{-1} + \frac{\eta\gamma}{2} \langle \theta_k^\eta, h(\theta_k^\eta) - H(\theta_k^\eta, X_{k+1})\rangle + \eta \langle h(\theta_k^\eta) - H(\theta_k^\eta,X_{k+1}), V_k^\eta \rangle
\\&  + \left( \frac{L_1 C_\rho \eta^2}{2} - \frac{\gamma^2 \lambda \eta^2}{4} - \frac{\gamma\eta}{2} + \frac{\gamma^2\eta^2}{4} + \frac{\gamma \lambda \eta}{2} \right) |V_k^\eta|^2  + \frac{\eta^2}{2} | H(\theta_k^\eta, X_{k+1})|^2 + \eta \gamma \beta^{-1} | \xi_{k+1} |^2 + \Sigma_k, \\
&\leq (1 - \gamma\lambda\eta) \frac{\cV_k}{\beta} + \eta \gamma A_c \beta^{-1} + \frac{\eta\gamma}{2} \langle \theta_k^\eta, h(\theta_k^\eta) - H(\theta_k^\eta, X_{k+1})\rangle + \eta \langle h(\theta_k^\eta) - H(\theta_k^\eta,X_{k+1}), V_k^\eta \rangle
\\&  + \left( \frac{L_1 C_\rho \eta^2}{2} - \frac{\gamma^2 \lambda \eta^2}{4} - \frac{\gamma\eta}{2} + \frac{\gamma^2\eta^2}{4} + \frac{\gamma \lambda \eta}{2} \right)  |V_k^\eta|^2  \\& + \frac{\eta^2}{2} \left( 2 L_1^2 (1 + |X_{k+1}|)^{2\rho} |\theta_k^\eta|^2 + 4 L_2^2 (1 + |X_{k+1}|)^{2(\rho+1)} + 4 H_0^2  \right)
+ \eta \gamma \beta^{-1} | \xi_{k+1} |^2 + \Sigma_k,
\end{align*}
by Remark~\ref{rem:BoundsOnH}. Let $\phi = (1 - \gamma \lambda \eta)$ and let ${\mathcal{H}_k} = \mathcal{G}_k \vee \sigma(\xi_1,\ldots,\xi_k)$. By using $\lambda \leq 1/4$, we obtain
\begin{align*}
&\frac{\bE[ \cV^2_{k+1} | {\mathcal{H}_k}]}{\beta^2} \leq \phi^2 \frac{\cV_k^2}{\beta^2} + 2 \phi \frac{\cV_k}{\beta} \left( \eta \gamma A_c \beta^{-1} + \left( \frac{L_1 C_\rho \eta^2}{2} - \frac{\gamma^2 \lambda \eta^2}{4} + \frac{\gamma^2\eta^2}{4} \right) |V_k^\eta|^2 \right.
\\& \left. + \frac{\eta^2}{2} \left( 2 L_1^2 C_\rho |\theta_k^\eta|^2 + 4 L_2^2 C_\rho + 4 H_0^2  \right) + \eta \gamma \beta^{-1} {d} \vphantom{\left( \frac{L}{2} \right)} \right) + \bE\left[\left(\vphantom{\left( \frac{L}{2} \right)} \eta \gamma A_c \beta^{-1} + \frac{\eta\gamma}{2} \langle \theta_k^\eta, h(\theta_k^\eta) - H(\theta_k^\eta, X_{k+1})\rangle \right.\right. \\
&\left. \left. \vphantom{\Big|} + \eta \langle h(\theta_k^\eta) - H(\theta_k^\eta,X_{k+1}), V_k^\eta \rangle + \left( \frac{L_1 C_\rho \eta^2}{2} - \frac{\gamma^2 \lambda \eta^2}{4} - \frac{\gamma\eta}{2} + \frac{\gamma^2\eta^2}{4} + \frac{\gamma \lambda \eta}{2} \right)  |V_k^\eta|^2  \right. \right. \\
& \left. \left. + \frac{\eta^2}{2} \left( 2 L_1^2 (1 + |X_{k+1}|)^{2\rho} |\theta_k^\eta|^2 + 4 L_2^2 (1 + |X_{k+1}|)^{2(\rho+1)} + 4 H_0^2  \right) + \eta \gamma \beta^{-1} | \xi_{k+1} |^2 + \Sigma_k \vphantom{\left( \frac{L_1 C_\rho \eta^2}{2} - \frac{\gamma^2 \lambda \eta^2}{4} + \frac{\gamma^2\eta^2}{4} \right)} \right)^2 \Big| {\mathcal{H}_k}\right].
\end{align*}
We first note that
\begin{align}
\frac{\cV_k}{\beta} &\geq \max\left\lbrace \frac{1}{8} (1 - 2\lambda) \gamma^2 |\theta_k^\eta|^2, \frac{1}{4} (1 - 2\lambda) |V_k^\eta|^2 \right\rbrace, \nonumber \\
&\geq \frac{1}{16} (1 - 2\lambda) \gamma^2 |\theta_k^\eta|^2 + \frac{1}{8} (1 - 2\lambda) |V_k^\eta|^2, \label{eq:cVLowerBound}
\end{align}
since $\max\{x,y\} \geq (x+y) / 2$ for any $x,y > 0$. Then we first obtain
\begin{align*}
&\frac{\bE[ \cV^2_{k+1} | {\mathcal{H}_k}]}{\beta^2} \leq (\phi^2 + 2 \phi \tilde{K}_1 \eta^2 ) \frac{\cV_k^2}{\beta^2} + 2 \phi \frac{\cV_k}{\beta} \left( \eta \gamma A_c \beta^{-1} + 2 \eta^2 L_2 C_\rho + 2 \eta^2 H_0^2 + \eta \gamma \beta^{-1} d \right)
\\&+ \bE\left[\left( \eta \gamma A_c \beta^{-1} + \frac{\eta\gamma}{2} \langle \theta_k^\eta, h(\theta_k^\eta) - H(\theta_k^\eta, X_{k+1})\rangle \right. \right. \\
& \left.\left. + \eta \langle h(\theta_k^\eta) - H(\theta_k^\eta,X_{k+1}), V_k^\eta \rangle + \left( \frac{L_1 C_\rho \eta^2}{2} - \frac{\gamma^2 \lambda \eta^2}{4} - \frac{\gamma\eta}{2} + \frac{\gamma^2\eta^2}{4} + \frac{\gamma \lambda \eta}{2} \right)  |V_k^\eta|^2  \right.\right. \\
& \left. \left. + \frac{\eta^2}{2} \left( 2 L_1^2 (1 + |X_{k+1}|)^{2\rho} |\theta_k^\eta|^2 + 4 L_2^2 (1 + |X_{k+1}|)^{2(\rho+1)} + 4 H_0^2  \right) + \eta \gamma \beta^{-1} | \xi_{k+1} |^2 + \Sigma_k \vphantom{\left( \frac{L_1 C_\rho \eta^2}{2} - \frac{\gamma^2 \lambda \eta^2}{4} + \frac{\gamma^2\eta^2}{4} \right)} \right)^2 \Big\vert {\mathcal{H}_k} \right]
\end{align*}
where
\begin{align*} 
\tilde{K}_1 := \max\left\lbrace \frac{\frac{{L}_1 C_\rho}{2} + \frac{ \gamma^2}{4} - \frac{\gamma^2 \lambda}{4}}{\frac{1}{8} (1 - 2\lambda)},  \frac{L_1^2 C_\rho}{\frac{1}{16} (1 - 2\lambda) \gamma^2} \right\rbrace
\end{align*}
Next using $\langle a,b \rangle \leq (|a|^2 + |b|^2)/2$, we obtain
\begin{align*}
&\frac{\bE[ \cV^2_{k+1} | {\mathcal{H}_k}]}{\beta^2} \leq (\phi^2 + 2 \phi \tilde{K}_1 \eta^2 ) \frac{\cV_k^2}{\beta^2} + 2 \phi \frac{\cV_k}{\beta} \left( \eta \gamma A_c \beta^{-1} + 2 \eta^2 L_2 C_\rho + 2 \eta^2 H_0^2 + \eta \gamma \beta^{-1} d \right)
\\&+ \bE\left[\left( \eta \gamma A_c \beta^{-1} + \frac{\eta\gamma}{4} |\theta_k^\eta|^2 + \left(\frac{\eta\gamma}{4} + \frac{\eta}{2}\right) | h(\theta_k^\eta) - H(\theta_k^\eta, X_{k+1}) |^2 \right. \right. \\
& \left.\left. + \left( \frac{L_1 C_\rho \eta^2}{2} - \frac{\gamma^2 \lambda \eta^2}{4} - \frac{\gamma\eta}{2} + \frac{\gamma^2\eta^2}{4} + \frac{\gamma \lambda \eta}{2} + \frac{\eta}{2} \right)  |V_k^\eta|^2  \right.\right. \\
& \left. \left. + \frac{\eta^2}{2} \left( 2 L_1^2 (1 + |X_{k+1}|)^{2\rho} |\theta_k^\eta|^2 + 4 L_2^2 (1 + |X_{k+1}|)^{2(\rho+1)} + 4 H_0^2  \right) + \eta \gamma \beta^{-1} | \xi_{k+1} |^2 + \Sigma_k \vphantom{\left( \frac{L_1 C_\rho \eta^2}{2} - \frac{\gamma^2 \lambda \eta^2}{4} + \frac{\gamma^2\eta^2}{4} \right)} \right)^2 \Big\vert {\mathcal{H}_k} \right]
\\&
\leq (\phi^2 + 2 \phi \tilde{K}_1 \eta^2 ) \frac{\cV_k^2}{\beta^2} + 2 \phi \frac{\cV_k}{\beta} \left( \eta \gamma A_c \beta^{-1} + 2 \eta^2 L_2 C_\rho + 2 \eta^2 H_0^2 + \eta \gamma \beta^{-1} d \right)
\\&+ \bE\left[\left( \eta \gamma A_c \beta^{-1} + \frac{\eta\gamma}{4} |\theta_k^\eta|^2 + \left(\frac{\eta\gamma}{2} + {\eta}\right) (| h(\theta_k^\eta)|^2 + | H(\theta_k^\eta, X_{k+1}) |^2) \right. \right. \\
& \left.\left. + \left( \frac{L_1 C_\rho \eta^2}{2} - \frac{\gamma^2 \lambda \eta^2}{4} - \frac{\gamma\eta}{2} + \frac{\gamma^2\eta^2}{4} + \frac{\gamma \lambda \eta}{2} + \frac{\eta}{2} \right)  |V_k^\eta|^2  \right.\right. \\
& \left. \left. + \frac{\eta^2}{2} \left( 2 L_1^2 (1 + |X_{k+1}|)^{2\rho} |\theta_k^\eta|^2 + 4 L_2^2 (1 + |X_{k+1}|)^{2(\rho+1)} + 4 H_0^2  \right) + \eta \gamma \beta^{-1} | \xi_{k+1} |^2 + \Sigma_k \vphantom{\left( \frac{L_1 C_\rho \eta^2}{2} - \frac{\gamma^2 \lambda \eta^2}{4} + \frac{\gamma^2\eta^2}{4} \right)} \right)^2 \Big\vert {\mathcal{H}_k} \right]
\end{align*}
where we have used $|a + b|^2 \leq 2|a|^2 + 2|b|^2$. Now using Remark~\ref{rem:BoundsOnH}, we obtain
\begin{align*}
&\frac{\bE[ \cV^2_{k+1} | {\mathcal{H}_k}]}{\beta^2} \leq (\phi^2 + 2 \phi \tilde{K}_1 \eta^2 ) \frac{\cV_k^2}{\beta^2} + 2 \phi \frac{\cV_k}{\beta} \left( \eta \gamma A_c \beta^{-1} + 2 \eta^2 L_2 C_\rho + 2 \eta^2 H_0^2 + \eta \gamma \beta^{-1} d \right)
\\&+ \bE\left[\left( \eta \gamma A_c \beta^{-1} + \frac{\eta\gamma}{4} |\theta_k^\eta|^2 + \left(\frac{\eta\gamma}{2} + {\eta}\right) \left( 2 L_1^2 C_{\rho}^2 |\theta_k^\eta|^2 + 2 h_0^2 + 2 L_1^2 (1 + |X_{k+1}|)^{2\rho} |\theta_k^\eta|^2 \right. \right. \right.
\\& \left.\left.\left. + 4 L_2^2 (1 + |X_{k+1}|)^{2(\rho+1)} + 4 H_0^2 \right) + \left( \frac{L_1 C_\rho \eta^2}{2} - \frac{\gamma^2 \lambda \eta^2}{4} - \frac{\gamma\eta}{2} + \frac{\gamma^2\eta^2}{4} + \frac{\gamma \lambda \eta}{2} + \frac{\eta}{2} \right)  |V_k^\eta|^2  \right.\right. \\
& \left. \left. + \frac{\eta^2}{2} \left( 2 L_1^2 (1 + |X_{k+1}|)^{2\rho} |\theta_k^\eta|^2 + 4 L_2^2 (1 + |X_{k+1}|)^{2(\rho+1)} + 4 H_0^2  \right) + \eta \gamma \beta^{-1} | \xi_{k+1} |^2 + \Sigma_k \vphantom{\left( \frac{L_1 C_\rho \eta^2}{2} - \frac{\gamma^2 \lambda \eta^2}{4} + \frac{\gamma^2\eta^2}{4} \right)} \right)^2 \Big\vert {\mathcal{H}_k} \right] \\
&=  (\phi^2 + 2 \phi \tilde{K}_1 \eta^2 ) \frac{\cV_k^2}{\beta^2} + 2 \phi \frac{\cV_k}{\beta} \left( \eta \gamma A_c \beta^{-1} + 2 \eta^2 L_2 C_\rho + 2 \eta^2 H_0^2 + \eta \gamma \beta^{-1} d \right)
\\&+ \bE\left[\left( \eta \gamma A_c \beta^{-1} + \left(\frac{\eta\gamma}{4} + \left(\frac{\eta\gamma}{2} + {\eta}\right) \left(2 L_1^2 C_{\rho}^2 +2 L_1^2 (1 + |X_{k+1}|)^{2\rho} \right) + \frac{\eta^2}{2} 2 L_1^2 (1 + |X_{k+1}|)^{2\rho} \right) |\theta_k^\eta|^2 \right. \right.
\\&
\left.\left. + \left(\frac{\eta\gamma}{2} + {\eta}\right) \left( 2 h_0^2 + 4 L_2^2 (1 + |X_{k+1}|)^{2(\rho+1)} + 4 H_0^2 \right) + \left( \frac{L_1 C_\rho \eta^2}{2} - \frac{\gamma^2 \lambda \eta^2}{4} - \frac{\gamma\eta}{2} + \frac{\gamma^2\eta^2}{4} + \frac{\gamma \lambda \eta}{2} + \frac{\eta}{2} \right)  |V_k^\eta|^2  \right.\right. \\
& \left. \left. + \frac{\eta^2}{2} \left( 4 L_2^2 (1 + |X_{k+1}|)^{2(\rho+1)} + 4 H_0^2  \right) + \eta \gamma \beta^{-1} | \xi_{k+1} |^2 + \Sigma_k \vphantom{\left( \frac{L_1 C_\rho \eta^2}{2} - \frac{\gamma^2 \lambda \eta^2}{4} + \frac{\gamma^2\eta^2}{4} \right)} \right)^2 \Big\vert {\mathcal{H}_k} \right]
\end{align*}
Now taking the expectation, regrouping and simplifying the terms above, we finally arrive at
\begin{align}
\frac{\bE[ \cV^2_{k+1} | {\mathcal{H}_k}]}{\beta^2} &\leq (\phi^2 + 2 \phi \tilde{K}_1 \eta^2 ) \frac{\cV_k^2}{\beta^2} + 2 \phi \frac{\cV_k}{\beta^2} \eta \tilde{c}_1 + 2 \phi \frac{\cV_k}{\beta} \eta^2 \hat{c}_1 + \eta^2 \frac{\tilde{c}_2}{\beta^2} + \tilde{c}_3 \eta^2 |\theta_k^\eta|^4 + \hat{c}_3 \eta^4 |\theta_k^\eta|^4 \nonumber \\ &
+ \eta^2 \tilde{c}_4 |V_k^\eta|^4 + \eta^4 \hat{c}_4 |V_k^\eta|^4 + \eta^2 \tilde{c}_5 + \eta^4 \hat{c}_5  + \eta^2 \frac{\tilde{c}_6}{\beta^2} + \bE[\Sigma_k^2 | {\mathcal{H}_k}]. \label{eq:IntermediateConditionalExp}
\end{align}
where
\begin{align*}
\tilde{c}_1 &= \gamma A_c + \gamma d, \\
\hat{c}_1 &= 2 L_2 C_\rho + 2 H_0^2 \\
\tilde{c}_2 &= 6 \gamma^2 A_c^2 \\
\tilde{c}_3 &= \frac{9\gamma^2}{8} + 18 (\gamma+2)^2 (L_1^4 C_{\rho}^4 + L_1^4 C_{\rho}), \\
\hat{c}_3 &= 18 L_1^4 C_{\rho}, \\
\tilde{c}_4 &= 6  \frac{(1 + \lambda \gamma - \gamma)^2}{4}, \\
\hat{c}_4 &= 6 \left( \frac{L_1 C_\rho}{2} + \frac{\gamma}{4} + \frac{1}{2}\right)^2, \\
\tilde{c}_5 &= (\gamma+2)^2 \left( 30 h_0^4 + 120 L_2^4 + 120 H_0^4 + 30 L_2^4 C_\rho + 30 H_0^4 \right), \\
\hat{c}_5 &= 48 (L_4^4 + H_0^4), \\
\tilde{c}_6 &= \gamma^2 d (d + 2).
\end{align*}
Next, recall
\begin{align*}
\Sigma_k = \frac{\gamma^2}{2} \sqrt{2\gamma^{-1}\beta^{-1} \eta} \langle \Delta_k^2, \xi_{k+1}\rangle + \frac{1}{2} \sqrt{2 \eta \gamma \beta^{-1}} \langle \Delta_k^1, \xi_{k+1}\rangle,
\end{align*}
where
\begin{align*}
\Delta_k^1 &= (1-\eta\gamma)V_k^\eta~-\eta~H(\theta_k^\eta, X_{k+1}) \\
\Delta_k^2 &= \theta_k^\eta + \gamma^{-1} V_k^\eta - \eta \gamma^{-1} H(\theta_k^\eta, X_{k+1}).
\end{align*}
Note that
\begin{align*}
\Sigma_k^2 &\leq 2 \gamma \beta^{-1} \eta |\Delta_k^1|^2 | \xi_{k+1}|^2 + 2 \gamma^3 \beta^{-1} \eta | \Delta_k^2|^2 |\xi_{k+1}|^2, \\
&\leq 2 \eta \gamma \beta^{-1} |\xi_{k+1}|^2 \left( |V_k^\eta - \eta (\gamma V_k^\eta + H(\theta_k^\eta, X_{k+1})|^2 + \gamma^2 |\theta_k^\eta + \gamma^{-1} V_k^\eta - \eta \gamma^{-1} H(\theta_k^\eta, X_{k+1})|^2 \right), \\
&\leq 2\eta\gamma\beta^{-1} |\xi_{k+1}|^2 \left( 2 (1 - \eta\gamma)^2 |V_k^\eta|^2 + 2 |H(\theta_k,X_{k+1})|^2 + 3\gamma^2 |\theta_k^\eta|^2 + 3 |V_k^\eta|^2 + 3 \eta^2 |H(\theta_k^\eta,X_{k+1})|^2 \right) \\
&\leq 2\eta\gamma\beta^{-1} |\xi_{k+1}|^2 \left( 2 (1 - \eta\gamma)^2 |V_k^\eta|^2 + 2 \left( 3 L_1^2 (1 + |X_{k+1}|)^{2\rho} |\theta_k^\eta|^2 + 3 L_2^2 (1 + |X_{k+1}|)^{2(\rho+1)} + 3 H_0^2 \right) \right. \\
&\left. + 3\gamma^2 |\theta_k^\eta|^2 + 3 |V_k^\eta|^2 + 3 \eta^2 \left(3 L_1^2 (1 + |X_{k+1}|)^{2\rho} |\theta_k^\eta|^2 + 3 L_2^2 (1 + |X_{k+1}|)^{2(\rho+1)} + 3 H_0^2\right) \right).
\end{align*}
Therefore, taking expectation and using \eqref{eq:cVLowerBound} and $\eta \leq 1/\gamma$, we obtain
\begin{align*}
\bE[\Sigma_k^2 | {\mathcal{H}_k}] &\leq 2\eta\gamma\beta^{-1} d \left( 5 |V_k^\eta|^2 + \left( 6 L_1^2 C_\rho + 3\gamma^2 + 9 L_1^2 C_\rho \right) |\theta_k^\eta|^2 + 15 L_2^2 C_\rho + 15 H_0^2 \right).
\\&
\leq  \eta \tilde{c}_7 \frac{\cV_k}{\beta^2} + \eta \frac{\tilde{c}_8}{\beta}
\end{align*}
where
\begin{align*} 
\tilde{c}_7 &= \max\left\lbrace \frac{10 \gamma d}{\frac{1}{8} (1 - 2\lambda)},  \frac{(2 \gamma d) (6 L_1^2 C_\rho + 3\gamma^2 + 9 L_1^2 C_\rho)}{\frac{1}{16} (1 - 2\lambda) \gamma^2} \right\rbrace \\
\tilde{c}_8 &= 30 \gamma d L_2^2 C_\rho + 30 \gamma d H_0^2.
\end{align*}
Now plugging this into the inequality \eqref{eq:IntermediateConditionalExp}, we arrive at
\begin{align}\label{eq:IntermediateCondExp2}
\frac{\bE[ \cV^2_{k+1} | {\mathcal{H}_k}]}{\beta^2} &\leq (\phi^2 + 2 \phi \tilde{K}_1 \eta^2 ) \frac{\cV_k^2}{\beta^2} + (2 \phi \tilde{c}_1 + \tilde{c}_7) \frac{\cV_k}{\beta^2} \eta+ 2 \phi \frac{\cV_k}{\beta} \eta^2 \hat{c}_1 + \eta^2 \frac{\tilde{c}_2}{\beta^2} + \tilde{c}_3 \eta^2 |\theta_k^\eta|^4 + \hat{c}_3 \eta^4 |\theta_k^\eta|^4 \nonumber \\
&+ \eta^2 \tilde{c}_4 |V_k^\eta|^4 + \eta^4 \hat{c}_4 |V_k^\eta|^4 + \eta^2 \tilde{c}_5 + \eta^4 \hat{c}_5  + \eta^2 \frac{\tilde{c}_6}{\beta^2} + \eta \frac{\tilde{c}_8}{\beta},
\end{align}
Now, recall our step-size condition for the above inequality, in particular, {$\eta \leq 1$}, which leads to
\begin{align}\label{eq:IntermediateCondExp3}
\frac{\bE[ \cV^2_{k+1} | {\mathcal{H}_k}]}{\beta^2} &\leq (\phi^2 + 2 \phi \tilde{K}_1 \eta^2 ) \frac{\cV_k^2}{\beta^2} + (2 \phi \tilde{c}_1 + \tilde{c}_7) \frac{\cV_k}{\beta^2} \eta + 2 \phi \frac{\cV_k}{\beta} \eta^2 \hat{c}_1 + \eta^2 \frac{\tilde{c}_2}{\beta^2} + (\tilde{c}_3 + \hat{c}_3) \eta^2 |\theta_k^\eta|^4 \nonumber \\
&+ \eta^2 (\tilde{c}_4 + \hat{c}_4) |V_k^\eta|^4 + \eta^2 (\tilde{c}_5 + \hat{c}_5) + \eta^2 \frac{\tilde{c}_6}{\beta^2} + \eta \frac{\tilde{c}_8}{\beta},
\end{align}
Next, recall that
\begin{align*}
\frac{\cV_k}{\beta} &\geq \max\left\lbrace \frac{1}{8} (1 - 2\lambda) \gamma^2 |\theta_k^\eta|^2, \frac{1}{4} (1 - 2\lambda) |V_k^\eta|^2 \right\rbrace,
\end{align*}
which implies
\begin{align}
\frac{\cV^2_k}{\beta^2} &\geq \max\left\lbrace \frac{1}{64} (1 - 2\lambda)^2 \gamma^4 |\theta_k^\eta|^4, \frac{1}{16} (1 - 2\lambda)^2 |V_k^\eta|^4 \right\rbrace, \nonumber \\
&\geq \frac{1}{128} (1 - 2\lambda)^2 \gamma^4 |\theta_k^\eta|^4 + \frac{1}{32} (1 - 2\lambda)^2 |V_k^\eta|^4. \label{eq:cV2LowerBound}
\end{align}
Now inequalities \eqref{eq:IntermediateCondExp3} and \eqref{eq:cV2LowerBound} together imply
\begin{align*}
\frac{\bE[ \cV^2_{k+1} | {\mathcal{H}_k}]}{\beta^2} &\leq (\phi^2 + 2 \phi \tilde{K}_2 \eta^2 ) \frac{\cV_k^2}{\beta^2} + (2 \phi \tilde{c}_1 + \tilde{c}_7) \frac{\cV_k}{\beta^2} \eta + 2 \phi \frac{\cV_k}{\beta} \eta^2 \hat{c}_1 + \eta^2 \frac{\tilde{c}_2}{\beta^2} + \eta^2 (\tilde{c}_5 + \hat{c}_5) + \eta^2 \frac{\tilde{c}_6}{\beta^2} + \eta \frac{\tilde{c}_8}{\beta},
\end{align*}
where $\tilde{K}_2 = \tilde{K_1} + \tilde{c}_{9}$ and
\begin{align*}
\tilde{c}_{9} = \max\left\lbrace \frac{\tilde{c}_3 + \hat{c}_3}{\frac{1}{128} (1 - 2\lambda)^2 \gamma^4},  \frac{\hat{c}_4 + \tilde{c}_4}{\frac{1}{32} (1 - 2\lambda)^2} \right\rbrace.
\end{align*}
Now we take unconditional expectations, use the fact that $\phi \leq 1$, and obtain
\begin{align*}
\bE[ \cV^2_{k+1}] &\leq \left(1 - {\lambda \gamma \eta} + \bar{K} \eta^2\right) \bE[\cV_k^2] + \tilde{c}_{10} \bE[\cV_k] \eta + 2 \hat{c}_1 \bE[\cV_k] \beta \eta^2 + \eta^2 \tilde{c}_2 + \eta^2 (\hat{c}_5 + \tilde{c}_5) \beta^2 + \eta^2 \tilde{c}_6 + \eta \tilde{c}_8 \beta,
\end{align*}
where $\bar{K} = 2 \tilde{K}_2$, $\tilde{c}_{10} = 2\tilde{c}_1 + \tilde{c}_7$. First note that, it follows from the proof of Lemma~\ref{lem:DiscreteTimeLemma} that
\begin{align*}
\bE[\cV_k] \leq \bE[\cV_0] +  \frac{4 (A_c + d)}{\gamma\lambda} := \tilde{c}_{11}.
\end{align*}
Therefore, using $\eta \leq 1$, we obtain
\begin{align*}
\bE[ \cV^2_{k+1}] &\leq \left(1 - {\lambda \gamma \eta} + \bar{K} \eta^2\right) \bE[\cV_k^2] + D \eta
\end{align*}
where $D = \tilde{c}_{10} \tilde{c}_{11} + 2 {\tilde{c}_{11}} \hat{c}_1 \beta + \tilde{c}_2 + (\hat{c}_5 + \tilde{c}_5) \beta^2 + \tilde{c}_6 + \tilde{c}_8 \beta$. Now recall that {$\eta \leq \lambda \gamma / 2 \bar{K}$} which leads to
\begin{align*}
\bE[ \cV^2_{k+1}] &\leq \left(1 - \frac{\lambda \gamma \eta}{2}\right) \bE[\cV_k^2] + D \eta,
\end{align*}
which implies that
\begin{align*}
\bE[\cV_k^2] \leq {\bE[\cV_0^2]}+ \frac{2 D}{\gamma\lambda},
\end{align*}
which concludes the proof. $\square$
\subsection{Proof of Lemma~\ref{lem:UpperLowerQuadBoundForU}}\label{proof:lem:UpperLowerQuadBoundForU}
We start by writing that
\begin{align*}
U(\theta) - U(0) &= \int_0^1 \langle \theta, h(t \theta) \rangle \md t, \\
&\leq \int_0^1 |\theta| |h(t\theta)| \md t, \\
&\leq \int_0^1 |\theta| {(\bar{L}_1 t |\theta| + h_0)} \md t.
\end{align*}
from Remark~\eqref{eq:BoundedGradh} where $\bar{L}_1 = L_1 \bE[(1+|X_0|)^\rho]$. This in turn leads to
\begin{align*}
U(\theta) &\leq u_0 + {\frac{\bar{L}_1}{2} |\theta|^2 + h_0 |\theta|.}
\end{align*}
where $u_0 = U(0)$. Next, we prove the lower bound. To this end, take $c \in (0,1)$ and write
\begin{align*}
U(\theta) &= U(c\theta) + \int_c^1 \langle \theta, h(t\theta)\rangle \md t, \\
&\geq \int_c^1 \frac{1}{t} \langle t\theta, h(t\theta) \rangle \md t, \\
&\geq \int_c^1 \frac{1}{t} \left(a |t\theta|^2 - b \right) \md t, \\
&= \frac{a (1 - c^2)}{2} |\theta|^2 + b \log c.
\end{align*}
Taking $c = 1 / \sqrt{3}$ leads to the bound.
\subsection{Proof of Lemma~\ref{lem:boundedVariance}}\label{proof:lem:boundedVariance}
Let $\mathcal{H}^\infty_t = \mathcal{F}^{\eta}_{\infty} \vee \mathcal{G}_{\lfloor t \rfloor}$. Following \cite{zhang2023nonasymptotic}, we obtain
\begin{align*}
&\E\left[\left|h(\bar{\zeta}_t^{\eta,n}) - H(\bar{\zeta}_t^{\eta,n},X_{nT+k})\right|^2\right] \\
& =\E\left[\E\left[\left.\left|h(\bar{\zeta}_t^{\eta,n}) - H(\bar{\zeta}_t^{\eta,n},X_{nT+k})\right|^2\right|\mathcal{H}^\infty_{nT}\right]\right] \\
&=\E\left[\E\left[\left.\left|\E\left[\left.H(\bar{\zeta}_t^{\eta,n}, X_{nT+k})\right|\mathcal{H}^\infty_{nT}\right] - H(\bar{\zeta}_t^{\eta,n},X_{nT+k})\right|^2\right|\mathcal{H}^\infty_{nT}\right]\right] \\
&\leq 4\E\left[\E\left[\left.\left| H(\bar{\zeta}_t^{\eta,n},X_{nT+k})- H(\bar{\zeta}_t^{\eta,n}, \E\left[\left. X_{nT+k}\right|\mathcal{H}_{nT}\right])\right|^2\right|\mathcal{H}_{nT}^\infty \right]\right] \\
&\leq 4L_2^2{\sigma}_Z\E\left[\left(1+\left|\bar{\zeta}_t^{\eta,n} \right|\right)^2\right],
\end{align*}
where the first inequality holds due to Lemma \ref{lem:magicLemma} and
\begin{align*}
\sigma_Z = \bE[(1 + |X_0| + | \bE[X_0]|)^{2\rho} |X_0 - \bE[X_0]|^2].
\end{align*}
Then, by using Lemma \ref{lem:zetaprocessmoment}, we obtain
\[
\E\left[\left|h(\bar{\zeta}_t^{\eta,n}) - H(\bar{\zeta}_t^{\eta,n},X_{nT+k})\right|^2\right] \leq
8 L_2^2 \sigma_Z + 8 L_2^2 \sigma_Z C_\zeta.
\]
\subsection{Proof of Lemma~\ref{lem:distanceToAuxiliary}}\label{proof:lem:distanceToAuxiliary}
Note that for any $t$, we have
\begin{align*}
\overline{V}_t^\eta = \overline{V}_\floor{t}^\eta - \eta\gamma \int_{\floor{t}}^t \overline{V}_{\floor{s}}^\eta \md s - \eta \int_{\floor{t}}^t H(\overline{\theta}_{\floor{s}},X_{\ceil{s}}) \md s + \sqrt{2 \gamma \eta \beta^{-1}} (B^\eta_t - B^\eta_{\floor{t}}).
\end{align*}
We therefore obtain
\begin{align*}
\bE\left[\left| \overline{V}_\floor{t}^\eta - \overline{V}_t^\eta \right|^2\right] &= \bE\left[\left| - \eta\gamma \int_{\floor{t}}^t \overline{V}_{\floor{s}}^\eta \md s - \eta \int_{\floor{t}}^t H(\overline{\theta}_{\floor{s}},X_{\ceil{s}}) \md s + \sqrt{2 \gamma \eta \beta^{-1}} (B^\eta_t - B^\eta_{\floor{t}}) \right|^2\right] \\
&\leq 4 \eta^2 \gamma^2 C_v + 4 \eta^2 (\widetilde{L}_1 C_\theta + \widetilde{C}_1) + 4 \gamma \eta \beta^{-1} d.
\end{align*}
Next, we write
\begin{align*}
\overline{\theta}_t^\eta = \overline{\theta}_{\floor{t}}^\eta + \int_{\floor{t}}^t \eta \overline{V}^\eta_{\floor{\tau}} \md \tau,
\end{align*}
which implies
\begin{align*}
\bE\left[\left|\overline{\theta}_t^\eta - \overline{\theta}_{\floor{t}}^\eta\right|^2\right] \leq \eta^2 C_v.
\end{align*}
\subsection{Proof of Lemma~\ref{lem:LowerBoundOfInnerProd}}\label{proof:lem:LowerBoundOfInnerProd}
Choose $\lambda > 0$ such that
\begin{align*}
\lambda = \min\left\lbrace \frac{1}{4}, \frac{a}{\bar{L}_1 + 2 \bar{L}_1 b + \frac{\gamma^2}{2}}\right\rbrace.
\end{align*}
Using Assumption~\ref{assum:dissipativity}, we obtain
\begin{align*}
\langle \theta, h(\theta) \rangle &\geq a |\theta|^2 - b, \\
&\geq 2\lambda \left( \frac{\bar{L}_1}{2} + \bar{L}_1 b + \frac{\gamma^2}{4} \right) |\theta|^2 - b, \\
&\geq 2 \lambda \left( U(\theta) - u_0 - \bar{L}_1 b |\theta| + \bar{L}_1 b |\theta|^2 + \frac{\gamma^2}{4} |\theta|^2 \right) - b, \\
&\geq 2 \lambda \left( U(\theta) - u_0 - \bar{L}_1 b + \frac{\gamma^2}{4} |\theta|^2 \right) - b,
\end{align*}
where the third line follows from Lemma~\ref{lem:UpperLowerQuadBoundForU} and the last line follows from the inequality $|x| \leq 1 + |x|^2$. Consequently, we obtain
\begin{align*}
\langle \theta, h(\theta) \rangle &\geq 2 \lambda \left( U(\theta) + \frac{\gamma^2}{4} |\theta|^2 \right) - 2 A_c/\beta
\end{align*}
where
\begin{align*}
A_c = \frac{\beta}{2} (b + 2 \lambda u_0 + 2 \lambda \bar{L}_1 b),
\end{align*}
which proves the claim.
\section{Proofs of main results}\label{app:sec:proofMainResults}

\subsection{Proof of Theorem~\ref{thm:W2_first_term}}\label{proof:thm:W2_first_term}
We note that
\begin{align}\label{eq:WassersteinDecomp}
W_2(\cL(\overline{\theta}_t^\eta,\overline{V}_t^\eta), \cL(\overline{\zeta}_t^{\eta,n},\overline{Z}_t^{\eta,n})) \leq \bE\left[\left|\overline{\theta}_t^\eta - \overline{\zeta}_t^{\eta,n}\right|^2\right]^{1/2} + \bE\left[\left|\overline{V}_t^{\eta} - \overline{Z}_t^{\eta,n}\right|^2\right]^{1/2}.
\end{align}
We first bound the first term of \eqref{eq:WassersteinDecomp}. We start by employing the synchronous coupling and obtain
\begin{align*}
\left| \overline{\theta}_t^{\eta} -\overline{\zeta}_t^{\eta,n} \right| &\leq \eta \int_{nT}^t \left| \overline{V}_\floor{s}^\eta -  \overline{Z}_s^{\eta,n} \right| \md s,
\end{align*}
which implies
\begin{align*}
\sup_{nT\leq u \leq t} \bE\left[ \left| \overline{\theta}_u^{\eta} -\overline{\zeta}_u^{\eta,n} \right|^2 \right] &\leq \eta \sup_{nT\leq u \leq t} \int_{nT}^u \bE\left[ \left| \overline{V}_\floor{s}^\eta -  \overline{Z}_s^{\eta,n} \right|^2 \right] \md s, \\
&= \eta \int_{nT}^t \bE\left[ \left| \overline{V}_\floor{s}^\eta -  \overline{Z}_s^{\eta,n} \right|^2 \right] \md s,
\end{align*}
Next, we write for any $t\in[nT,(n+1)T)$
\begin{align*}
\left| \overline{V}_\floor{t}^\eta -  \overline{Z}_t^{\eta,n} \right|
&\leq \left| \overline{V}_{\floor{t}}^\eta - \overline{V}_t^\eta \right| +  \left| \overline{V}_t^\eta - \overline{Z}_t^{\eta,n} \right|, \\
& \leq \left| \overline{V}_{\floor{t}}^\eta - \overline{V}_t^\eta \right| + \left| -\gamma\eta \int_{nT}^t [\overline{V}_\floor{s}^\eta - \overline{Z}_s^{\eta,n}] \md s - \eta \int_{nT}^t [H(\overline{\theta}_\floor{s}^\eta, X_{\ceil{s}}) - h(\overline{\zeta}_s^{\eta,n})] \md s \right|, \\
&\leq \left| \overline{V}_{\floor{t}}^\eta - \overline{V}_t^\eta \right| + \gamma\eta \int_{nT}^t \left| \overline{V}_\floor{s}^\eta - \overline{Z}_s^{\eta,n} \right| \md s + \eta \left| \int_{nT}^t [H(\overline{\theta}_\floor{s}^\eta, X_{\ceil{s}}) - h(\overline{\zeta}_s^{\eta,n})] \md s \right|, \\
&\leq \left| \overline{V}_{\floor{t}}^\eta - \overline{V}_t^\eta \right| + \gamma\eta \int_{nT}^t \left| \overline{V}_\floor{s}^\eta - \overline{Z}_s^{\eta,n} \right| \md s \\ &+ \eta \int_{nT}^t \left| H(\overline{\theta}_\floor{s}^\eta, X_{\ceil{s}}) - H(\overline{\zeta}_s^{\eta,n},X_\ceil{s})\right| \md s + \eta \left| \int_{nT}^t [H(\overline{\zeta}_s^{\eta,n}, X_{\ceil{s}}) - h(\overline{\zeta}_s^{\eta,n})] \md s \right|, \\
&\leq \left| \overline{V}_{\floor{t}}^\eta - \overline{V}_t^\eta \right| + \gamma\eta \int_{nT}^t \left| \overline{V}_\floor{s}^\eta - \overline{Z}_s^{\eta,n} \right| \md s \\
&+ \eta \int_{nT}^t L_1 (1 + X_\ceil{s})^\rho \left| \overline{\theta}_\floor{s}^\eta - \overline{\zeta}_s^{\eta,n} \right| \md s + \eta \left| \int_{nT}^t [H(\overline{\zeta}_s^{\eta,n}, X_{\ceil{s}}) - h(\overline{\zeta}_s^{\eta,n})] \md s \right|.
\end{align*}
We take the squares of both sides and use $(a+b)^2 \leq {2 (a^2 +  b^2)}$ twice to obtain
\begin{align*}
\left| \overline{V}_\floor{t}^\eta -  \overline{Z}_t^{\eta,n} \right|^2 &\leq 4 \left| \overline{V}_{\floor{t}}^\eta - \overline{V}_t^\eta \right|^2 + 4 \gamma^2 \eta^2 \left( \int_{nT}^t \left| \overline{V}_\floor{s}^\eta - \overline{Z}_s^{\eta,n} \right| \md s \right)^2 \\
&+ 4\eta^2 \left( \int_{nT}^t L_1 (1 + X_{\ceil{s}})^\rho \left| \overline{\theta}_{\floor{s}}^\eta - \overline{\zeta}_s^{\eta,n} \right| \md s \right)^2 + 4 \eta^2 \left| \int_{nT}^t [H(\overline{\zeta}_s^{\eta,n},X_\ceil{s}) - h(\overline{\zeta}_s^{\eta,n})] \md s \right|^2 \\
&\leq 4 \left| \overline{V}_{\floor{t}}^\eta - \overline{V}_t^\eta \right|^2 + 4 \gamma^2 \eta \int_{nT}^t \left| \overline{V}_\floor{s}^\eta - \overline{Z}_s^{\eta,n} \right|^2 \md s \\
&+ 4\eta L_1^2 \int_{nT}^t (1 + X_{\ceil{s}})^{2\rho} \left| \overline{\theta}_{\floor{s}}^\eta - \overline{\zeta}_s^{\eta,n} \right|^2 \md s + 4 \eta^2 \left| \int_{nT}^t [H(\overline{\zeta}_s^{\eta,n},X_\ceil{s}) - h(\overline{\zeta}_s^{\eta,n})] \md s \right|^2.
\end{align*}
Taking expectations of both sides, we obtain
\begin{align*}
\bE\left[ \left| \overline{V}_\floor{t}^\eta -  \overline{Z}_t^{\eta,n} \right|^2 \right]&\leq 4 \bE\left[ \left| \overline{V}_{\floor{t}}^\eta - \overline{V}_t^\eta \right|^2\right] + 4 \gamma^2 \eta \int_{nT}^t \bE\left[ \left| \overline{V}_\floor{s}^\eta - \overline{Z}_s^{\eta,n} \right|^2\right] \md s \\
&+ 4\eta L_1^2 C_\rho \int_{nT}^t \bE\left[\left| \overline{\theta}_{\floor{s}}^\eta - \overline{\zeta}_s^{\eta,n} \right|^2\right] \md s + 4 \eta^2 \bE\left[ \left| \int_{nT}^t [H(\overline{\zeta}_s^{\eta,n},X_\ceil{s}) - h(\overline{\zeta}_s^{\eta,n})] \md s \right|^2 \right], \\
&\leq 4 {\sigma}_V \eta + 4 \gamma^2 \eta \int_{nT}^t \bE\left[ \left| \overline{V}_\floor{s}^\eta - \overline{Z}_s^{\eta,n} \right|^2\right] \md s \\
&+ 4\eta L_1^2 C_\rho \int_{nT}^t \bE\left[\left| \overline{\theta}_{\floor{s}}^\eta - \overline{\zeta}_s^{\eta,n} \right|^2\right] \md s + 4 \eta^2 \bE\left[ \left| \int_{nT}^t [H(\overline{\zeta}_s^{\eta,n},X_\ceil{s}) - h(\overline{\zeta}_s^{\eta,n})] \md s \right|^2 \right].
\end{align*}
By applying Gr\"{o}nwall's lemma, we arrive at
\begin{align*}
\bE\left[ \left| \overline{V}_\floor{t}^\eta -  \overline{Z}_t^{\eta,n} \right|^2 \right] &\leq 4 c_1 {\sigma}_V \eta + 4\eta c_1 L_1^2 C_\rho \int_{nT}^t \bE\left[\left| \overline{\theta}_{\floor{s}}^\eta - \overline{\zeta}_s^{\eta,n} \right|^2\right] \md s  \\ & +4 c_1 \eta^2 \bE\left[ \left| \int_{nT}^t [H(\overline{\zeta}_s^{\eta,n},X_\ceil{s}) - h(\overline{\zeta}_s^{\eta,n})] \md s \right|^2 \right],
\end{align*}
where $c_1 = \exp( 4\gamma^2)$ since $\eta T \leq 1$. Next, we write
\begin{align}
\sup_{nT \leq u \leq t} \bE\left[ \left| \overline{\theta}_u^{\eta} -\overline{\zeta}_u^{\eta,n} \right|^2 \right] &\leq \eta \int_{nT}^t \bE\left[ \left| \overline{V}_\floor{s}^\eta -  \overline{Z}_s^{\eta,n} \right|^2 \right] \md s, \nonumber \\
&\leq 4 c_1 \eta {\sigma}_V + 4 {c_1} \eta^2 L_1^2 C_\rho \int_{nT}^t \int_{nT}^s \bE\left[\left| \overline{\theta}_{\floor{s'}}^\eta - \overline{\zeta}_{s'}^{\eta,n} \right|^2\right] \md s' \md s \nonumber \\
& + 4 c_1 \eta^3 \int_{nT}^t \bE\left[ \left| \int_{nT}^s [H(\overline{\zeta}_{s'}^{\eta,n},X_\ceil{s'}) - h(\overline{\zeta}_{s'}^{\eta,n})] \md s' \right|^2 \right] \md s, \nonumber \\
&\leq 4 c_1 \eta {\sigma}_V + 4 {c_1} \eta L_1^2 C_\rho \sup_{nT \leq s \leq t} \int_{nT}^s \bE\left[\left| \overline{\theta}_{\floor{s'}}^\eta - \overline{\zeta}_{s'}^{\eta,n} \right|^2\right] \md s' \nonumber \\
& + 4 c_1 \eta^3 \int_{nT}^t \bE\left[ \left| \int_{nT}^s [H(\overline{\zeta}_{s'}^{\eta,n},X_\ceil{s'}) - h(\overline{\zeta}_{s'}^{\eta,n})] \md s' \right|^2 \right] \md s. \label{eq:MainIntermediateIneq}
\end{align}
First, we bound the supremum term (i.e. the second term of \eqref{eq:MainIntermediateIneq}) as
\begin{align}
\sup_{nT \leq s \leq t} \int_{nT}^s \bE \left[ \left| \overline{\theta}_{\floor{s'}}^\eta - \overline{\zeta}_{s'}^{\eta,n}\right|^2\right] \md s' &= \int_{nT}^t \bE\left[\left| \overline{\theta}_{\floor{s'}}^\eta - \overline{\zeta}_{s'}^{\eta,n} \right|^2\right] \md s' \nonumber \\
&\leq \int_{nT}^t \sup_{nT \leq u \leq s'} \bE\left[\left| \overline{\theta}_\floor{u}^\eta - \overline{\zeta}_u^{\eta,n}\right|^2\right] \md s' \nonumber \\
&\leq \int_{nT}^t \sup_{nT \leq u \leq s'} \bE\left[\left| \overline{\theta}_u^\eta - \overline{\zeta}_u^{\eta,n}\right|^2\right] \md s'. \label{eq:SupremumBound}
\end{align}
Next, we bound the last term of \eqref{eq:MainIntermediateIneq} by partitioning the integral. Assume that $nT + K \leq s \leq t \leq nT + K + 1$ where $K + 1 \leq T$. Thus we can write
\begin{align*}
\left| \int_{nT}^s \left[ h(\bar{\zeta}_{s'}^{\eta,n}) - H(\bar{\zeta}_{s'}^{\eta,n},X_{\ceil{{s'}}})\right] \md {s'} \right| = \left| \sum_{k=1}^{K} I_k + R_K \right|
\end{align*}
where
\begin{align*}
I_k = \int_{nT + (k-1)}^{nT + k} [h(\bar{\zeta}_{s'}^{\eta,n}) - H(\bar{\zeta}_{s'}^{\eta,n},X_{nT+k})] \md {s'} \quad \quad \textnormal{and} \quad\quad R_K = \int_{nT+K}^{s} [h(\bar{\zeta}_{s'}^{\eta,n}) - H(\bar{\zeta}_{s'}^{\eta,n},X_{nT+K + 1})] \md {s'}.
\end{align*}
Taking squares of both sides
\begin{align*}
\left| \sum_{k=1}^{K} I_k + R_K \right|^2 = \sum_{k=1}^K | I_k |^2 + 2 \sum_{k=2}^K \sum_{j = 1}^{k-1} \langle I_k, I_j \rangle +2 \sum_{k=1}^K \langle I_k, R_K \rangle + |R_K|^2,
\end{align*}
Finally, it remains to take the expectations of both sides. We begin by defining the filtration $\mathcal{H}^\infty_s = \mathcal{F}^{\eta}_{\infty} \vee \mathcal{G}_{\lfloor s \rfloor}$ and note that for any $k =2, \dots, K$, $j = 1, \dots, k-1$,
\begin{align*}
&\bE \langle I_k, I_j \rangle\\
 &= \bE\left[ \bE [\langle I_k, I_j \rangle | \mathcal{H}^\infty_{nT+j} ] \right], \\
&= \bE\left[ \bE \left[\left\langle \int_{nT + (k-1)}^{nT + k} [H(\bar{\zeta}_{s'}^{\eta,n},X_{nT + k}) - h(\bar{\zeta}_{s'}^{\eta,n})] \md {s'}, \left.\int_{nT + (j-1)}^{nT + j} [H(\bar{\zeta}_{s'}^{\eta,n},X_{nT + j}) - h(\bar{\zeta}_{s'}^{\eta,n})] \md {s'} \right\rangle \right|  \mathcal{H}^\infty_{nT+j} \right] \right], \\
& = \bE\left[ \left\langle \int_{nT + (k-1)}^{nT + k} \bE \left[\left. H(\bar{\zeta}_{s'}^{\eta,n},X_{nT + k}) - h(\bar{\zeta}_{s'}^{\eta,n})\right|  \mathcal{H}^\infty_{nT+j} \right]\md {s'},  \int_{nT + (j-1)}^{nT + j} [H(\bar{\zeta}_{s'}^{\eta,n},X_{nT + j}) - h(\bar{\zeta}_{s'}^{\eta,n})] \md {s'} \right\rangle  \right], \\
&= 0.
\end{align*}
By the same argument $\bE \langle I_k, R_K\rangle = 0$ for all $1 \leq k \leq K$. Therefore,
\begin{align}
\int_{nT}^t &\bE\left[ \left| \int_{nT}^s [H(\overline{\zeta}_{s'}^{\eta,n},X_\ceil{s'}) - h(\overline{\zeta}_{s'}^{\eta,n})] \md s' \right|^2 \right] \md s \nonumber\\
&= \int_{nT}^t \left[ \sum_{k=1}^K \bE \left[\left| \int_{nT + (k-1)}^{nT + k} [h(\bar{\zeta}_{s'}^{\eta,n}) - H(\bar{\zeta}_{s'}^{\eta,n},X_{nT+k})] \md {s'}\right|^2 \right] \right] \md s \nonumber\\
&+\int_{nT}^t \bE\left[ \left| \int_{nT+K}^{s} [h(\bar{\zeta}_{s'}^{\eta,n}) - H(\bar{\zeta}_{s'}^{\eta,n},X_{nT+K + 1})] \md {s'} \right|^2 \right] \md s \nonumber\\
& \leq \int_{nT}^t \left[ \sum_{k=1}^K \int_{nT + (k-1)}^{nT + k} \bE\left[ \left| h(\bar{\zeta}_{s'}^{\eta,n}) - H(\bar{\zeta}_{s'}^{\eta,n},X_{nT+k})\right|^2 \right] \md {s'} \right] \md s \nonumber\\
&+\int_{nT}^t \int_{nT+K}^{s} \bE\left[ \left| h(\bar{\zeta}_{s'}^{\eta,n}) - H(\bar{\zeta}_{s'}^{\eta,n},X_{nT+K + 1})\right|^2\right] \md {s'} \md s \nonumber\\
& \leq T^2 \sigma_H + T \sigma_H. \label{eq:IntegralsSplitBound}
\end{align}
Using \eqref{eq:MainIntermediateIneq}, \eqref{eq:SupremumBound}, and \eqref{eq:IntegralsSplitBound}, we eventually obtain
\begin{align}
\sup_{nT \leq u \leq t} \bE\left[ \left| \overline{\theta}_u^{\eta} -\overline{\zeta}_u^{\eta,n} \right|^2 \right] &\leq 4 c_1 \eta {\sigma}_V + 4 \eta L_1^2 C_\rho \int_{nT}^t \sup_{nT \leq u \leq s'} \bE\left[\left| \overline{\theta}_u^\eta - \overline{\zeta}_u^{\eta,n}\right|^2\right] \md s', \nonumber \\
& + 4 c_1 \eta^3 (T^2 \sigma_H + T \sigma_H), \nonumber \\
&\leq 4 c_1 \eta {\sigma}_V + 4 {c_1} \eta L_1^2 C_\rho \int_{nT}^t \sup_{nT \leq u \leq s'} \bE\left[\left| \overline{\theta}_u^\eta - \overline{\zeta}_u^{\eta,n}\right|^2\right] \md s' \nonumber \\ & + 4 c_1 \eta \sigma_H + 4 c_1 \eta^2 \sigma_H, \label{eq:MainIntermediateIneq2}
\end{align}
since $\eta T \leq 1$. Finally, applying Gr\"onwall's inequality and using again $\eta T \leq 1$ provides
\begin{align*}
\sup_{nT \leq u \leq t} \bE\left[ \left| \overline{\theta}_u^{\eta} -\overline{\zeta}_u^{\eta,n} \right|^2 \right] \leq \exp(4 {c_1}L_1^2 C_\rho) (4 c_1 \sigma_V + 4 c_1 \sigma_H + 4 c_1 \eta \sigma_H) \eta,
\end{align*}
which implies that
\begin{align}\label{eq:WassersteinDecompBoundFirst}
\bE\left[ \left| \overline{\theta}_t^{\eta} -\overline{\zeta}_t^{\eta,n} \right|^2 \right]^{1/2} \leq C_{1,1}^\star \sqrt{\eta}
\end{align}
with $C_{1,1}^\star = \sqrt{\exp(4 {c_1} L_1^2 C_\rho) (4 c_1 \sigma_V + 4 c_1 \sigma_H + 4 c_1 \eta \sigma_H)}$. Note that $\sigma_V = \mathcal{O}(d)$ and $\sigma_H = \mathcal{O}(d)$ hence $C_{1,1}^\star = \mathcal{O}(\sqrt{d})$.

Next, we upper bound the second term of \eqref{eq:WassersteinDecomp}. To prove it, we write
\begin{align*}
\left| \overline{V}_t^\eta - \overline{Z}_t^{\eta,n}\right| \leq \left| \gamma \eta \int_{nT}^t \left[ \overline{V}_{\floor{s}}^\eta - \overline{Z}_s^{\eta,n}\right] \md s \right| + \eta \left| \int_{nT}^t \left[H(\overline{\theta}^\eta_{\floor{s}}) - h(\overline{\zeta}_s^{\eta,n})\right] \md s \right|,
\end{align*}
which leads to
\begin{align*}
\bE\left[ \left| \overline{V}_t^\eta - \overline{Z}_t^{\eta,n}\right|^2\right] \leq 2 \gamma^2 \eta \int_{nT}^t \bE\left[ \left| \overline{V}_{\floor{s}}^\eta - \overline{Z}_s^{\eta,n}\right|^2 \right] + 2 \eta^2 \bE\left[ \left| \int_{nT}^t \left[H(\overline{\theta}^\eta_{\floor{s}}) - h(\overline{\zeta}_s^{\eta,n})\right] \md s \right|^2\right].
\end{align*}
By similar arguments we have used for bounding the the first term, we obtain
\begin{align*}
\bE\left[ \left| \overline{V}_t^\eta - \overline{Z}_t^{\eta,n}\right|^2\right] \leq 2 \gamma^2 \eta \int_{nT}^t \bE\left[ \left| \overline{V}_{\floor{s}}^\eta - \overline{Z}_s^{\eta,n}\right|^2 \right] + {\exp(4 c_1 L_1^2 C_\rho) (4 c_1 \sigma_H \eta + 4 c_1 \eta^2 \sigma_H \eta^2).}
\end{align*}
Using the fact that the rhs is an increasing function of $t$ and we obtain
\begin{align*}
\sup_{nT \leq u \leq t} \bE\left[ \left| \overline{V}_u^\eta - \overline{Z}_u^{\eta,n}\right|^2\right] \leq 2 \gamma^2 \eta \int_{nT}^t \sup_{nT\leq u \leq s} \bE\left[ \left| \overline{V}_{u}^\eta - \overline{Z}_u^{\eta,n}\right|^2 \right] \md s + {\exp(4 c_1 L_1^2 C_\rho) (4 c_1 \sigma_H \eta + 4 c_1 \eta^2 \sigma_H \eta^2)},
\end{align*}
Applying Gronwall's lemma and $\eta T \leq 1$ yields
\begin{align*}
\sup_{nT \leq u \leq t} \bE\left[ \left| \overline{V}_u^\eta - \overline{Z}_u^{\eta,n}\right|^2\right] \leq  \exp(2\gamma^2 + {4c_1 L_1^2 C_\rho}) 4 c_1 (\sigma_H \eta + \sigma_H \eta^2),
\end{align*}
which leads to
\begin{align}\label{eq:WassersteinDecompSecondTerm}
\bE\left[ \left| \overline{V}_t^\eta - \overline{Z}_t^{\eta,n}\right|^2\right]^{1/2} \leq C_{1,2}^\star \sqrt{\eta},
\end{align}
where $C_{1,2}^\star = \sqrt{\exp(2\gamma^2 + {4c_1 L_1^2 C_\rho}) 4 c_1 \sigma_H ( 1 + \eta)}$. Note again that $\sigma_H = \mathcal{O}(d)$, hence $C_{1,2}^\star = \mathcal{O}(d^{1/2})$.

Therefore, combining \eqref{eq:WassersteinDecomp}, \eqref{eq:WassersteinDecompBoundFirst}, \eqref{eq:WassersteinDecompSecondTerm}, we obtain
\begin{align*}
W_2(\cL(\overline{\theta}_t^\eta,\overline{V}_t^\eta), \cL(\overline{\zeta}_t^{\eta,n},\overline{Z}_t^{\eta,n})) \leq C_1^\star \eta^{1/2},
\end{align*}
where $C_1^\star = C_{1,1}^\star + C_{1,2}^\star = \mathcal{O}(d^{1/2})$.
\subsection{Proof of Theorem~\ref{thm:W2_second_term}}\label{proof:thm:W2_second_term}
{Triangle inequality implies that
\begin{align*}
W_2&(\cL(\overline{\zeta}_t^{\eta,n},\overline{Z}_t^{\eta,n}), \cL(\zeta_t^\eta,Z_t^\eta)) \leq \sum_{k=1}^n W_2(\cL(\overline{\zeta}_t^{\eta,k},\overline{Z}_t^{\eta,k}), \cL(\overline{\zeta}_t^{\eta,k-1},\overline{Z}_t^{\eta,k-1})), 
\end{align*}
since $(\overline{\zeta}_t^{\eta,0}, \overline{Z}_t^{\eta,0}) = (\widehat{\zeta}_t^{0, \overline{\theta}_0^\eta, \overline{V}_0^\eta, \eta}, \widehat{Z}_t^{0, \overline{\theta}_0^\eta, \overline{V}_0^\eta, \eta}) = (\zeta_t^\eta, Z_t^\eta)$ by definition. Next we write out the definition of the process $(\overline{\zeta}_t^{\eta,k},\overline{Z}_t^{\eta,k})$ and obtain
\begin{align*}
W_2&(\cL(\overline{\zeta}_t^{\eta,n},\overline{Z}_t^{\eta,n}), \cL(\zeta_t^\eta,Z_t^\eta)) \\ &\leq \sum_{k=1}^n W_2(\cL(\widehat{\zeta}_t^{kT,\overline{\theta}_{kT}^\eta,\overline{V}_{kT}^\eta,\eta},\widehat{Z}_t^{kT,\overline{\theta}_{kT}^\eta,\overline{V}_{kT}^\eta,\eta}), \cL(\widehat{\zeta}_t^{(k-1)T,\overline{\theta}_{(k-1)T}^\eta,\overline{V}_{(k-1)T}^\eta,\eta},\widehat{Z}_t^{(k-1)T,\overline{\theta}_{(k-1)T}^\eta,\overline{V}_{(k-1)T}^\eta,\eta})).
\end{align*}
At this point, we have two processes and in order to be able to use the contraction result, we need their starting times to match. For notational simplicity, let us define
\begin{align*}
\widehat{B}_t^{s, \overline{\theta}^\eta_s, \overline{V}^\eta_s, \eta} = (\widehat{\zeta}_t^{s,\overline{\theta}_{s}^\eta,\overline{V}_{s}^\eta,\eta},\widehat{Z}_t^{s,\overline{\theta}_{s}^\eta,\overline{V}_{s}^\eta,\eta})
\end{align*}
Therefore, in order to be able to use a contraction result, we note
\begin{align*}
\cL(\widehat{B}_t^{
(k-1)T, \overline{\theta}^\eta_{(k-1)T}, \overline{V}^\eta_{(k-1)T}, \eta}) = 
\cL(\widehat{B}_t^{
kT, \widehat{B}_{kT}^{(k-1)T, \overline{\theta}^\eta_{(k-1)T}, \overline{V}^\eta_{(k-1)T}, \eta}, \eta}).
\end{align*}
This leads to
\begin{align*}
&W_2(\cL(\overline{\zeta}_t^{\eta,n},\overline{Z}_t^{\eta,n}), \cL(\zeta_t^\eta,Z_t^\eta)) \\ &\leq \sum_{k=1}^n  W_2(\cL(\widehat{\zeta}_t^{kT,\theta_{kT}^\eta,V_{kT}^\eta,\eta},\widehat{Z}_t^{kT,\theta_{kT}^\eta,V_{kT}^\eta,\eta}),\cL(\widehat{\zeta}_t^{kT,\widehat{B}_{kT}^{(k-1)T, \overline{\theta}^\eta_{(k-1)T}, \overline{V}^\eta_{(k-1)T}, \eta},\eta},\widehat{Z}_t^{kT,\widehat{B}_{kT}^{(k-1)T, \overline{\theta}^\eta_{(k-1)T}, \overline{V}^\eta_{(k-1)T}, \eta},\eta})).
\end{align*}
The reader should notice at this point that this quantity can be upper bounded by a \textit{contraction result} as both processes are defined for time $t$ and both started at time $kT$. By using Theorem~\ref{thm:Eberle}, we obtain
\begin{align*}
&W_2(\cL(\overline{\zeta}_t^{\eta,n},\overline{Z}_t^{\eta,n}), \cL(\zeta_t^\eta,Z_t^\eta)) \leq \sqrt{\dot{C}} \sum_{k=1}^n e^{-\eta\dot{c}(t-kT)/2} \sqrt{\mathcal{W}_\rho(\cL(\theta_{kT}^\eta,V_{kT}^\eta),\cL(\overline{\zeta}_{kT}^{\eta,(k-1)T}, \overline{Z}_{kT}^{\eta, (k-1)T}}).
\end{align*}
Next, using Lemma~5.4 of \cite{chau2022stochastic} to upper bound the the last term leads
\begin{align*}
&W_2(\cL(\overline{\zeta}_t^{\eta,n},\overline{Z}_t^{\eta,n}), \cL(\zeta_t^\eta,Z_t^\eta)) \\ 
&\leq 3 \max\{1 + \alpha, \gamma^{-1}\} \sqrt{\dot{C}} \sum_{k=1}^n e^{-\eta\dot{c}(t-kT)/2} \sqrt{1 + \varepsilon_c \bE^{1/2}[\mathcal{V}^2(\theta_{kT}^\eta,V_{kT}^\eta)] + \varepsilon_c \bE^{1/2}[\mathcal{V}^2(\overline{\zeta}_{kT}^{\eta,(k-1)T},\overline{Z}_{kT}^{\eta,(k-1)T})]} \\
&\times \sqrt{W_2(\cL(\theta_{kT}^\eta,V_{kT}^\eta),\cL(\overline{\zeta}_{kT}^{\eta,(k-1)T},\overline{Z}_{kT}^{\eta,(k-1)T}))}, \\
&\leq C_2^\star \eta^{1/4},
\end{align*}
where the last line follows from Lemma~\ref{lem:SquareContraction} and Theorem~\ref{thm:W2_first_term}. The bound for the term relating to the continuous dynamics (i.e. $\bE^{1/2}[\mathcal{V}^2(\overline{\zeta}_{kT}^{\eta,(k-1)T},\overline{Z}_{kT}^{\eta,(k-1)T})]$) follows from Remark~\ref{R2} (dissipativity of the gradient) and we have a constant diffusion coefficient, and uniform bound in Lemma~\ref{lem:ContTimeLemma} and, finally, the fact that the initial data (which is the iterates of the numerical scheme) has finite fourth moments (Lemma~\ref{lem:SquareContraction}). We note that $C_2^\star = \mathcal{O}(e^{d})$.}
\subsection{Proof of Proposition~\ref{prop:BoundT1}}\label{proof:prop:BoundT1}
We denote $\pi_{n,\beta}^\eta := \cL(\theta_n^\eta,V_n^\eta)$ and write
\begin{align*}
\bE[U(\theta_n^\eta)] - \bE[U(\theta_\infty)] = \int_{\bR^{2d}} U(\theta) \pi_{n,\beta}^\eta(\md \theta,\md v) - \int_{\bR^{2d}} U(\theta) \pi_\beta(\md \theta, \md v).
\end{align*}
Recall from \eqref{eq:BoundedGradh}, we have
\begin{align*}
|h(\theta)| \leq \overline{L}_1 | \theta| + h_0.
\end{align*}
{Using \citet[Lemma~6]{raginsky2017non}, we arrive at}
\begin{align*}
\left| \int_{\bR^{2d}} U(\theta) \pi_{n,\beta}^\eta(\md \theta,\md v) - \int_{\bR^{2d}} U(\theta) \pi_\beta(\md \theta, \md v) \right| \leq (\overline{L}_1 {C_m} + h_0) W_2(\pi_{n,\beta}^\eta,\pi_\beta),
\end{align*}
where,
{\begin{align*}
C^2_m := \max\left( \int_{\bR^{2d}} \|\theta\|^2 \pi_{n,\beta}^\eta(\md \theta, \md v), \int_{\bR^{2d}} \|\theta\|^2 \pi_\beta(\md \theta, \md v)  \right) = \max(C_\theta^c, C_\theta).
\end{align*}
We therefore obtain using Theorem~\ref{thm:ConvRate} that
\begin{align*}
\bE[U(\theta_n^\eta)] - \bE[U(\theta_\infty)] \leq (\overline{L}_1 {C_m} + h_0) \left(C_1^\star d^{1/2} \eta^{1/2} + C_2^\star \eta^{1/4} + C_3^\star e^{-C_4^\star \eta n}\right).
\end{align*}}

\end{document}